\documentclass{amsart}
\usepackage{graphicx}
\usepackage{calc}
\usepackage{hyperref}

\usepackage[cmtip,arrow,all,2cell]{xy}
\usepackage{pb-diagram,pb-xy}
\usepackage{amsmath,amsthm,amssymb,tabularx}

\newcommand{\abs}[1]{\left\vert#1\right\vert}

\theoremstyle{plain}

\newtheorem{tm}{Theorem}[section]

\newtheorem{prop}[tm]{Proposition}
\newtheorem{lm}[tm]{Lemma}
\newtheorem{cor}[tm]{Corollary}

\theoremstyle{definition}

\newtheorem{ex}{Example}[section]

\newtheorem{rem}[ex]{Remark}

\newcommand{\btm}{\begin{tm}}
\newcommand{\etm}{\end{tm}}
\newcommand{\blm}{\begin{lm}}
\newcommand{\elm}{\end{lm}}
\newcommand{\bprop}{\begin{prop}}
\newcommand{\eprop}{\end{prop}}
\newcommand{\bcor}{\begin{cor}}
\newcommand{\ecor}{\end{cor}}
\newcommand{\bex}{\begin{ex}}
\newcommand{\eex}{\end{ex}}
\newcommand{\brem}{\begin{rem}}
\newcommand{\erem}{\end{rem}}

\newcommand{\bpr}{\begin{proof}}
\newcommand{\epr}{\end{proof}}
\newcommand{\beq}{\begin{equation}}
\newcommand{\eeq}{\end{equation}}
\newcommand{\bit}{\begin{itemize}}
\newcommand{\eit}{\end{itemize}}

\def\R{\mathbb{R}}
\def\C{\mathbb{C}}

\def\e{\varepsilon}

\def \le {\leqslant}

\def \ph {\varphi}
\def\id{\operatorname{id}}

\def\Ste{\operatorname{\sf Ste}}

\begin{document}

\title{On tensor fractions and tensor products in the category of stereotype spaces}

\author{S.S.Akbarov}
 \address{School of Applied Mathematics, National Research University Higher School of Economics, 34, Tallinskaya St. Moscow, 123458 Russia}
 \email{sergei.akbarov@gmail.com}

\thanks{Supported by the RFBR grant No. 18-01-00398.}

\maketitle

\begin{abstract}
We prove two identities that connect some natural tensor products in the category $\sf{LCS}$ of locally convex spaces with the tensor products in the category $\sf{Ste}$ of stereotype spaces. In particular, we give sufficient conditions under which the identity
$$
X^\vartriangle\odot Y^\vartriangle\cong (X^\vartriangle\cdot Y^\vartriangle)^\vartriangle\cong (X\cdot Y)^\vartriangle
$$
holds, where $\odot$ is the injective tensor product in the category $\sf{Ste}$,
$\cdot$, the primary tensor product in the category $\sf{LCS}$, and $\vartriangle$, the pseudosaturation operation in the category $\sf{LCS}$. Studying the relations of this type is justified by the fact that they turn out to be important instruments for constructing duality theory based on the notion of envelope. In particular, they are used in the construction of the duality theory for the class of (not necessarily, Abelian) countable discrete groups.
\end{abstract}

\section{Introduction}

In \cite{Ak95-1,Ak95-2,Ak03,Ak08,Ak16} the author described the properties of the class $\Ste$ of locally convex spaces named {\it stereotype spaces}, and defined by the condition\footnote{In G.~K\"othe's book \cite{Kothe} the spaces of this type are called {\it polar reflexive}.}
$$
X\cong (X^\star)^\star,
$$
where each star $\star$ means the dual space of functionals, endowed with the topology of uniform convergence on totally bounded sets. As was noticed, it has a series of remarkable properties, in particular,
\bit{

\item[---] $\Ste$ is very wide, since it contains all quasicomplete barreled spaces (in particular, all Fr\'echet spaces, and thus, all Banach spaces),

\item[---] $\Ste$ forms a pre-Abelian complete and co-complete category with linear continuous maps as morphisms,

\item[---] $\Ste$ is an autodual category with respect to the functor $X\mapsto X^\star$,

\item[---] $\Ste$ has three natural bifunctors $\circledast$, $\odot$, $\oslash$ with the following properties:
\begin{eqnarray}
 \C\circledast X \cong X  \cong  X\circledast \C, & \C\odot X \cong  X  \cong  X\odot \C, \label{C-circledast-X-cong-X}
\\
 X\circledast Y \cong  Y\circledast X, & X\odot Y \cong Y\odot
X, \label{X-circledast-Y-cong-Y-circledast-X}
\\
 (X\circledast Y)\circledast Z  \cong  X\circledast
(Y\circledast Z), &  (X\odot Y)\odot Z  \cong  X\odot (Y\odot Z),
\label{X-circledast-Y-circledast-Z}
 \\
    Z\oslash (Y\circledast X)\cong (Z\oslash Y)\oslash X,
 &  (X\odot Y)\oslash Z \cong X\odot (Y\oslash Z), \label{Z-oslash-(Y-circledast-X)} \\
(X\circledast Y)^\star  \cong   Y^\star \odot X^\star, &
(X\odot Y)^\star \cong   Y^\star \circledast X^\star, \\
 X\circledast Y\cong (X^\star\oslash Y)^\star, & X\odot Y\cong Y\oslash (X^\star).
\end{eqnarray}
}\eit
Identities \eqref{C-circledast-X-cong-X}---\eqref{X-circledast-Y-circledast-Z} in this list mean that $\Ste$ is a symmetric monoidal category with respect to bifunctors $\circledast$ and $\odot$, and this justifies the name  {\it ``stereotype tensor products''} for them ($\circledast$, the {\it projective tensor product}, and $\odot$, the {\it injective}). On the other hand, identities \eqref{Z-oslash-(Y-circledast-X)} justify the name {\it ``stereotype tensor fraction''} for the bifunctor $\oslash$ (and the left one means that with respect to $\circledast$ the monoidal category $\Ste$ is closed with the inner hom-functor $\oslash$).

In the further studies the author (and later, Yu.~N.~Kuznetsova and O.~Yu.~Aristov) found numerous applications of stereotype spaces in functional analysis and in geometry, in particular, in generalizations of the Pontryagin duality on non-commutative groups. As an example, from the results of \cite{Ak08,ArHRC} it follows that the Pontryagin duality is extended from the category of finite Abelian groups to the category of affine algebraic groups by the functor $G\mapsto {\mathcal O}^\star(G)$ of taking group algebra of analytic functionals. This generalization can be expressed by the following functorial diagram:
{\sf
 \beq\label{diagramma-kategorij-dlya-konechnyh-affinnyh-grupp}
 \xymatrix @R=3.pc @C=2.pc
 {
 \boxed{\begin{matrix}
  \text{Hopf algebras,}\\
  \text{reflexive with respect}\\
  \text{to the Arens---Michael envelope}
 \end{matrix}}
 \ar[rr]^{H\mapsto H^\dagger} & &
 \boxed{\begin{matrix}
  \text{Hopf algebras,}\\
  \text{reflexive with respect}\\
  \text{to the Arens---Michael envelope}
 \end{matrix}}
 \\
 \boxed{\begin{matrix}
  \text{affine algebraic groups}
 \end{matrix}} \ar[u]^(0.4){\scriptsize\begin{matrix} {\mathcal O}^\star(G)\\
 \text{\rotatebox{90}{$\mapsto$}} \\ G\end{matrix}} & &
 \boxed{\begin{matrix}
  \text{affine algebraic groups}
 \end{matrix}} \ar[u]_(0.4){\scriptsize\begin{matrix} {\mathcal O}^\star(G) \\
 \text{\rotatebox{90}{$\mapsto$}} \\ G\end{matrix}} \\
 \boxed{\begin{matrix}
 \text{Abelian finite groups}
 \end{matrix}} \ar[u]^{\mathfrak{e}} \ar[rr]^{G\mapsto G^\bullet} & &
  \boxed{\begin{matrix}
 \text{Abelian finite groups}
  \end{matrix}}\ar[u]_{\mathfrak{e}}
 }
 \eeq }\noindent
where $\mathfrak{e}$ means the natural embedding of categories, $G^\bullet$, the Pontryagin dual group, $\dagger$, the composition of the functor $H\mapsto \widehat{H}$ of the Arens---Michael envelope (see definition in \cite{Ak16} or in \cite{AkCountGroups}) and the functor $X\mapsto X^\star$ of taking dual stereotype space,
$$
H^\dagger:=( \widehat{H} )^\star,
$$
and by Hopf algebras reflexive with respect to the Arens---Michael envelope, we mean the Hopf algebra in the monoidal category $(\Ste,\circledast)$ satisfying the identity
$$
H\cong \big(H^\dagger\big)^\dagger
$$
(see details in \cite{AkCountGroups}). The commutativity of the diagram \eqref{diagramma-kategorij-dlya-konechnyh-affinnyh-grupp} is understood as an isomorphism of two functors goping from the left lower corner to the right upper one:
\beq\label{O^star(G)^dagger-cong-O^star(G^bullet)}
{\mathcal O}^\star(G)^\dagger\cong {\mathcal O}^\star(G^\bullet).
\eeq

After the author's work \cite{Ak08}, where the duality with respect to an envelope was first described and applied to complex Lie groups, the same idea (with the replacement of the Arens---Michael envelope by other envelopes) was used to the generalizations of the Pontryagin dualities to various other classes of groups. In particular, Yu.~N.~Kuznetsova in \cite{Kuz} suggested a duality theory for (not necessarily commutative) Moore groups (later her results were specified in \cite{Ak17-1,Ak17-2,Akbarov-Moore}), and, besides this, the author himself in \cite{Ak17-1,Ak17-2} described the duality theory for differential geometry, that works at least for the real Lie groups of the form $\R^n\times K\times D$, where $K$ is a compact Lie group, and $D$ a discrete Moore group.

These observations have not yet taken shape in general theories with natural boundaries, since it is obvious that all these results can be generalized to broader classes of groups, and it is not clear where it would be logical to stop in these studies. Within the framework of this program, in the work \cite{AkCountGroups}, the author considered the question of whether it is possible, by analogy with diagram \eqref{diagramma-kategorij-dlya-konechnyh-affinnyh-grupp} to generalize the Pontryagin duality from the category of finite Abelian groups to the category of countable discrete groups. In \cite{AkCountGroups} it was answered positively, with the describing diagram
{\sf
 \beq\label{diagramma-kategorij-dlya-diskretnyh-grupp}
 \xymatrix @R=3.pc @C=2.pc
 {
 \boxed{ \begin{matrix}
  \text{holomorphically reflexive}\\
  \text{Hopf algebras}\\
 \end{matrix}}
 \ar[rr]^{H\mapsto H^\dagger} & &
 \boxed{ \begin{matrix}
  \text{holomorphically reflexive}\\
  \text{Hopf algebras}\\
 \end{matrix}}
 \\
 \boxed{\begin{matrix}
 \text{countable discrete groups}
 \end{matrix}} \ar[u]^(0.4){\scriptsize\begin{matrix} {\mathcal O}^\star(G)\\
 \text{\rotatebox{90}{$\mapsto$}} \\ G\end{matrix}} & &
 \boxed{\begin{matrix}
 \text{countable discrete groups}
 \end{matrix}} \ar[u]_(0.4){\scriptsize\begin{matrix} {\mathcal O}^\star(G) \\
 \text{\rotatebox{90}{$\mapsto$}} \\ G\end{matrix}} \\
 \boxed{\begin{matrix}
 \text{Abelian finite groups}
 \end{matrix}} \ar[u]^{\mathfrak{e}} \ar[rr]^{G\mapsto G^\bullet} & &
  \boxed{\begin{matrix}
 \text{Abelian finite groups}
  \end{matrix}}\ar[u]_{\mathfrak{e}}
 }
 \eeq }\noindent
where the holomorphic reflexivity of Hopf algebras is understood as reflexivity with respect to one of the stereotype variants of the Arens---Michael envelope (see details in \cite{AkCountGroups}).

The proof of this fact, however, as it turns out, requires a deep immersion into the foundations of the theory of stereotype spaces with the establishing a series of technical statements, which, due to their volume and ideological isolation, should be presented as a separate text. This article is such a text. Here we derive two formulas expressing the injective stereotype tensor product $\odot$ and the tensor fraction $\oslash$, mentioned above, in terms of some auxiliary constructions in the theory of topological vector spaces.

The fist formula expresses $\odot$ in terms of one of the variants of tensor producst in the category of locally convex spaces, which we denote by the dot $\cdot$ and call a ``primary tensor product''. Formally $X\cdot Y$ is defined as the space of linear continuous mappings $\ph:X^\star\to Y$ with the topology of uniform convergence on polars $U^\circ$ of neighbourhoods of zero in $X$, and for it we prove the formula
\beq\label{(X^vartriangle-cdot-Y^vartriangle)^vartriangle=(X-cdot-Y)^-vartriangle-0}
X^\vartriangle\odot Y^\vartriangle\cong (X^\vartriangle\cdot Y^\vartriangle)^\vartriangle\cong (X\cdot Y)^\vartriangle.
\eeq
which holds for all pseudocomplete locally convex spaces $X$ and $Y$ (without the assumption that one of them possesses the classical approximation property).

And the second formula completes the formal definition of the tensor fraction $\oslash$ by the spaces of operators. Recall that in \cite{Ak03} the space $\oslash$ was defined in two steps. First, it was defined the space $Y:X$ (the ``primary tensor product''), which was the space of linear continuous mappings $\ph:X\to Y$ with the topology of uniform convergence on totally bounded sets. After that the space $Y\oslash X$ is defined as the pseudosaturation of the space $Y:X$
\beq\label{DEF:oslash}
    Y\oslash X=(Y:X)^\vartriangle.
\eeq
And here is another formula connecting $\oslash$ with $:$, and proved in our article:
\beq\label{(vartriangle_Y:triangledown_X)^vartriangle-isomorphism-0}
Y^\vartriangle\oslash X^\triangledown\cong(Y^\vartriangle:X^\triangledown)^\vartriangle\cong  (Y:X)^\vartriangle
\eeq
($\triangledown$ is the operation of pseudocompletion in the category of locally convex spaces \cite[1.3]{Ak03}, and $\vartriangle$ is again the operation of pseudosaturation).

The importance of formula \eqref{(X^vartriangle-cdot-Y^vartriangle)^vartriangle=(X-cdot-Y)^-vartriangle-0} is that it allows to to study the space $X^\vartriangle\odot Y^\vartriangle$ without the detailed study of the pseudosaturations $X^\vartriangle$ and $Y^\vartriangle$. For instance, in \cite{AkCountGroups} they are used in the proof of identities like
\beq\label{F(M)-odot-F(N)}
{\mathcal F}(M)^\vartriangle\odot {\mathcal F}(N)^\vartriangle\cong {\mathcal F}(M\times N)^\vartriangle
\eeq
where $M$ and $N$ are manifolds of certain type, and ${\mathcal F}(...)$ a certain type of (topological vector) functional spaces on these manifolds. A typical situation in the theory of stereotype spaces is when we easily can describe the topology of a given space, in this case, of ${\mathcal F}(...)$, but the topology of its pseudosaturation ${\mathcal F}(...)^\vartriangle$ looks extremely complicated, so that for its description it is usually impossible to find anything simpler than the very definition of the pseudosaturation. And formula \eqref{F(M)-odot-F(N)} gives the description of stereotype tensor product without the necessity to analyze the topologies of the spaces ${\mathcal F}(...)^\vartriangle$.

Similarly, formula \eqref{(vartriangle_Y:triangledown_X)^vartriangle-isomorphism-0} allows to reduce the description of the space $Y^\vartriangle\oslash X^\triangledown$ to the description of the space $Y:X$ without the necessity to study the spaces $Y^\vartriangle$ and $X^\triangledown$ (however, we give formula \eqref{(vartriangle_Y:triangledown_X)^vartriangle-isomorphism-0} for the future, since directly it is not used in \cite{AkCountGroups}).

These are the questions the present work is devoted to.

\subsection{Acknowledgements}

The author thanks O.~Yu.~Aristov for useful discussions.

\section{Pseudosaturation of the primary fraction $Y:X$}

Everywhere in the text we use the terminology and the notations of \cite{Ak03}. We consider locally convex spaces over the field $\C$. For each locally convex space $X$ the symbol ${\mathcal U}(X)$ denotes the set of all neighbourhoods of zero in $X$, ${\mathcal S}(X)$, the set of all totally bounded subsets in $X$, ${\mathcal D}(X)$ the set of all capacious\footnote{A set $D\subseteq X$ is said to be {\it capatious}, if for each totally bounded set $S\subseteq X$ there is a finite set $F\subseteq X$ such that $S\subseteq D+F$.} sets in $X$. The term {\it operator} will be used for linear continuous mappings $\ph:X\to Y$ of locally convex spaces.

We start with the proof of formula \eqref{(vartriangle_Y:triangledown_X)^vartriangle-isomorphism-0}, the last in the list of Introduction.

\subsection{The bilinear form $(g,x)\mapsto g:x$.}

Recall that in \cite[5.6]{Ak03} a bilinear mapping $\beta:X\times Y\to Z$ was said to be {\it continuous}, if\label{DEF:neprer-bilin-forma}
\bit{
\item[1)] for each neighbourhood of zero $W\subseteq Z$ and for exach totally bounded set $S\subseteq X$ there is a neighbourhood of zero $V\subseteq Y$ such that
    $$
    \beta(S,V)\subseteq W,
    $$
    and
\item[2)] for each neighbourhood of zero $W\subseteq Z$ and for exach totally bounded set $T\subseteq Y$ there is a neighbourhood of zero $U\subseteq X$ such that
    $$
    \beta(U,T)\subseteq W.
    $$

}\eit

Recall again that for each locally convex spaces $X$ and $Y$ the symbol $Y:X$ ({\it ``primary tensor fraction''}) denotes the space of linear continuous mappings $\ph:X\to Y$ with the topology of {\it uniform convergence on totally bounded sets} \cite[5.1]{Ak03}. If $A\subseteq X$ and $B\subseteq Y$, then the symbol $B:A$ denotes the subset in $Y:X$ consisting of mappings that turn $A$ into $B$ \cite[5.4]{Ak03}:
\beq\label{DEF:B:A}
\ph\in B:A\quad\Leftrightarrow\quad \ph\in Y:X\ \& \ \ph(A)\subseteq B.
\eeq
For each $x\in X$ and $g\in Y^\star$ we define the functional
\beq\label{DEF:(g:x)}
    (g:x):(Y:X)\to\C\quad\Big|\quad (g:x)(\ph)=g(\ph(x)),\qquad \ph\in Y:X
\eeq

\btm\label{TH:(g,x)|->g:x-razdelno-nepreryvno}
Suppose $X$ and $Y$ are locally convex spaces, and $X$ is pseudosaturated. Then the bilinear mapping
$$
(g,x)\in Y^\star\times X\mapsto g:x\in (Y:X)^\star
$$
is continuous (in the above sense).
\etm
\bpr
1. Take $g_i\to 0$ and $S\in{\mathcal S}(X)$. Then for arbitrary totally bounded set $\varPhi\subseteq Y:X$ the set $\varPhi(S)$ is totally bounded in $Y$, by \cite[Theorem 5.1]{Ak03}. So the net $g_i$ tends to zero uniformly on the set $\varPhi(S)$ (by definition of topology in $Y^\star$):
$$
(g_i:x)(\ph)=g_i(\ph(x)) \overset{\C}{\underset{i\to\infty}{\underset{\ph\in\varPhi,\ x\in S}{\rightrightarrows}}}0.
$$
This is true for each totally bounded set $\varPhi\subseteq Y:X$, hence
$$
g_i:x \overset{Y:X}{\underset{i\to\infty}{\underset{x\in S}{\rightrightarrows}}}0.
$$

2. Take $G\in {\mathcal S}(Y^\star)$ and $x_i\to 0$. Then for arbitrary totally bounded set $\varPhi\subseteq Y:X$ the set $G\circ\varPhi$ is totally bounded in $X^\star$ by \cite[Theorem 5.1]{Ak03}. Since $X$ is pseudosaturated, the set $G\circ\varPhi$ is equicontinuous on $X$. As a corollary,
$$
(g:x_i)(\ph)=g(\ph(x_i)) \overset{\C}{\underset{i\to\infty}{\underset{\ph\in\varPhi,\ g\in G}{\rightrightarrows}}}0.
$$
This is true for each totally bouded set $\varPhi\subseteq Y:X$, hence
$$
g:x_i \overset{Y:X}{\underset{i\to\infty}{\underset{g\in G}{\rightrightarrows}}}0.
$$
\epr

\bcor\label{COR:(g,x)|->g:x-razdelno-nepreryvno}
Suppose $X$ and $Y$ are locally convex spaces, and $X$ is pseudosaturated. If $S\subseteq X$ is totally bounded and $g\in Y^\star$, then  $g:S\subseteq (Y:X)^\star$ is totally bounded. Dually, if $x\in X$ and $G\subseteq Y^\star$ is totally bounded, then $G:x\subseteq (Y:X)^\star$ is totally bounded.
\ecor

\subsection{Pseudosaturation of $Y:X$.}

Let $\{X^\lambda;\ \lambda\in\bf{Ord}\}$ be an injective series of the space $X$, and  $\{Y_\mu;\ \mu\in\bf{Ord}\}$ a projective series of the space $Y$, and
$$
\vee_\lambda^X: X\to X^\lambda,\qquad \wedge_\mu^Y: Y_\mu\to Y
$$
the corresponding natural mappings. Consider the analog of the diagram \cite[(5.1)]{Ak03}:
\beq\label{wedge_mu^Y:vee_lambda^X}
\xymatrix @R=2.pc @C=2.0pc 
{
X \ar[d]_{(\wedge_\mu^Y:\vee_\lambda^X)(\psi)= \wedge_\mu^Y\circ\psi\circ \vee_\lambda^X} \ar[r]^{\vee_\lambda^X} & X^\lambda \ar[d]^{\psi} \\
Y & Y_\mu \ar[l]_{\wedge_\mu^Y}
}
\eeq
It defines a mapping
\beq\label{wedge_mu^Y:vee_lambda^X-1}
(\wedge_\mu^Y:\vee_\lambda^X): (Y_\mu:X^\lambda) \to (Y:X).
\eeq

\blm\label{LM:Y_mu:X^lambda=Y:X}
If $X$ is pseudosaturated and $Y$ is pseudocomplete, then the mapping \eqref{wedge_mu^Y:vee_lambda^X-1} is a bijection, so the spaces $Y_\mu:X^\lambda$ and $Y:X$ can be identified as sets (while as locally convex spaces they can have different topologies):
\beq\label{Y_mu:X^lambda=Y:X}
Y_\mu:X^\lambda=Y:X
\eeq
\elm
\bpr
The mapping $\vee_\lambda^X: X\to X^\lambda$ is an epimorphism. On the other hand, the mapping $\wedge_\mu^Y: Y_\mu\to Y$ is an injection. This implies that the mapping
$$
\psi\mapsto \wedge_\mu^Y\circ\psi\circ (\wedge_\lambda^X)^\star
$$
is an injection. We have to verify that it is a surjection. Indeed, let $\ph:X\to Y$ be an orbitrary operator. We subsequently construct two operators $\ph_\mu$ and $\psi$ such that the following diagram is commutative:
$$
\xymatrix @R=3.pc @C=4.0pc 
{
X \ar[d]_{\ph}\ar@{-->}[dr]_{\ph_\mu} \ar[r]^{\vee_\lambda^X} & X^\lambda \ar@{-->}[d]^{\psi} \\
Y & Y_\mu \ar[l]^{\wedge_\mu^Y}
}
$$
Since $X$ is pseudosaturated, the operator $\ph$ is uniquely extended to the operator $\ph^\vartriangle:X\to Y^\vartriangle$. Its composition with the natural mapping  $Y^\vartriangle\to Y_\mu$ is an operator that extends $\ph$ to an operator $\ph_\mu:X\to Y_\mu$.

Further, since $Y$ is pseudocomplete, $Y_\mu$ is again pseudocomplete (by \cite[Proposition 3.16]{Ak03}). Hence, the operator $\ph_\mu:X\to Y_\mu$ is extended to an operator $\chi:X^\triangledown\to Y_\mu$. Its composition with the natural embedding $X^\lambda\to X^\triangledown$ is the operator $\psi:X^\lambda\to Y_\mu$.
\epr

\blm\label{LM:o-sohranenii-polnoi-ogranichennosti-v-Y_mu:X_lambda}
Let $X$ be pseudosaturated, and $Y$ pseudocomplete. If $\varPhi\in Y:X$ is totally bounded, then for any $\lambda,\mu\in\bf{Ord}$ its representation in the space $Y_\mu:X^\lambda$ (by the bijection \eqref{Y_mu:X^lambda=Y:X}) is again a totally bounded set.
\elm
\bpr
1. Consider first the case of $\lambda=0$. If $\varPhi\subseteq Y:X$ is totally bounded, then by \cite[Theorem 5.1]{Ak03} this means that $\varPhi$ is equicontinuous and uniformly totally bounded on each totally bounded set $S\subseteq X$. As a corollary, the image $\varPhi(S)$ is totally bounded in $Y$. The space $Y_\mu$ is a strengthening of the topology of $Y$, but the class of totally bounded sets in $Y$ and the topology on them are not changed. As a corollary, $\varPhi(S)$ is totally bounded in the space $Y_\mu$ as well, and moreover, with the same topology as the one induced from $Y$. We can conclude that the set $\varPhi$, being considered as the set of operators from $X$ into $Y_\mu$, is equicontinuous and uniformly totally bounded on $S$. And this is true for each totally bounded set $S\subseteq X$. Hence (again by \cite[Theorem 5.1]{Ak03}), $\varPhi$ is totally bounded in $Y_\mu:X$.

2. Thus, we understood that if $\varPhi\subseteq Y:X$ is totally bounded in $Y:X$, then it is totally bounded in $Y_\mu:X$ for each $\mu\in\bf{Ord}$. If we take a big enough ordinal $\mu\in\bf{Ord}$, then we obtain that $\varPhi\subseteq Y:X$ is totally bounded in the space $Y^\vartriangle:X$. From this by \cite[Lemma 5.10]{Ak03} we obtain that $\varPhi\subseteq Y:X$ is totally bounded in the space $Y^{\vartriangle\triangledown}:X^\triangledown$.

Since, by assumption, $Y$ is pseudocomplete, by \cite[Proposition 3.16]{Ak03}, the space $Y^{\vartriangle}$ is also pseudocomplete, and, as a corollary,
$$
Y^{\vartriangle\triangledown}=Y^{\vartriangle}.
$$
We obtain that the set $\varPhi\subseteq Y:X$ is totally bounded in the space $Y^\vartriangle:X^\triangledown$.

Now if we take arbitrary $\lambda,\mu\in{\bf Ord}$ and denote by
$$
\sigma:X^\lambda\to X^\triangledown,\qquad
\pi:Y_\mu\gets Y^\vartriangle,
$$
the natural mappings, we obtain a linear continuous mapping
$$
(\pi:\sigma): (Y^\vartriangle:X^\triangledown) \to (Y_\mu:X^\lambda).
$$
It turns the totally bounded set $\varPhi\subseteq Y^\vartriangle:X^\triangledown$ into the totally bounded set $\varPhi\subseteq Y_\mu:X^\lambda$.
\epr

\blm\label{LM:o-drobyah-v-Y_mu:X_lambda}
Let $X$ be pseudosaturated, and $Y$ pseudocomplete. If $S\in {\mathcal S}(X^\lambda)$ and $V\in {\mathcal U}(Y_\mu)$, then $V:S\in {\mathcal D}(Y:X)$.\footnote{The notation $V:S$ is introduced in \eqref{DEF:B:A}.}
\elm
\bpr
Let $\varPhi\subseteq Y:X$ be a totally bounded set. By Lemma  \ref{LM:o-sohranenii-polnoi-ogranichennosti-v-Y_mu:X_lambda}, $\varPhi$ is totally bounded in the space $Y_\mu:X^\lambda$ as well. Since $V:S$ is a neighbourhood of zero in the space $Y_\mu:X^\lambda$, there is a finite set $A\subseteq Y_\mu:X^\lambda$ such that \beq\label{o-drobyah-v-Y_mu:X_lambda}
\varPhi\subseteq V:S+A.
\eeq
By Lemma \ref{LM:Y_mu:X^lambda=Y:X}, all these sets can be considered as subsets in $Y:X$. We obtain the following: for each totally bounded set $\varPhi$ in $Y:X$ there is a finite set $A$ in $Y:X$ such that \eqref{o-drobyah-v-Y_mu:X_lambda} holds. This means that the set $V:S$ is capacious in the space $Y:X$.
\epr

\blm\label{LM:(Y_mu:X_lambda)^vartriangle=(Y:X)^vartriangle}
If $X$ is pseudosaturated and $Y$ is pseudocomplete, then for any ordinals  $\lambda,\mu\in\bf{Ord}$
\bit{
\item[(i)] the spaces $Y_\mu:X^\lambda$ and $Y:X$ become isomorphic after pseudosaturation:
\beq\label{(Y_mu:X_lambda)^vartriangle=(Y:X)^vartriangle}
(Y_\mu:X^\lambda)^\vartriangle=(Y:X)^\vartriangle,
\eeq

\item[(ii)] for each $x\in X^\lambda$ and $g\in (Y_\mu)^\star$ the functional
$$
    (g:x):(Y_\mu:X^\lambda)\to\C\quad\Big|\quad (g:x)(\ph)=g(\ph(x)),\qquad \ph\in Y_\mu:X^\lambda
$$
    is continuous on the space $(Y:X)^\vartriangle$, i.e. there is a unique functional $h\in ((Y:X)^\vartriangle)^\star$ such that the following diagram is commutative:
\beq\label{(Y_mu:X_lambda)^vartriangle=(Y:X)^vartriangle-func}
\xymatrix @R=2.pc @C=4.0pc 
{
(Y_\mu:X^\lambda)^\vartriangle \ar[d]_{\vartriangle_{Y_\mu:X^\lambda}} \ar[r]^{(\wedge_\mu^Y:\vee_\lambda^X)^\vartriangle} & (Y:X)^\vartriangle \ar@{-->}[d]^{h} \\
Y_\mu:X^\lambda \ar[r]_{g:x}& \C
}
\eeq
}\eit
\elm
\bpr
We have to organize a double induction by $\lambda$ and $\mu$. For $\lambda=\mu=0$ equality \eqref{(Y_mu:X_lambda)^vartriangle=(Y:X)^vartriangle-func} becomes trivial
$$
(Y:X)^\vartriangle=(Y:X)^\vartriangle,
$$
and diagram \eqref{(Y_mu:X_lambda)^vartriangle=(Y:X)^vartriangle-func} turns into the diagram
$$
\xymatrix @R=2.pc @C=8.0pc 
{
(Y:X)^\vartriangle \ar[d]_{\vartriangle_{Y:X}} \ar[r]^{(\id_Y:\id_X)^\vartriangle=\id_{Y:X}^\vartriangle} & (Y:X)^\vartriangle \ar@{-->}[d]^{h} \\
Y:X \ar[r]_{g:x}& \C
}
$$
where we just have to put $h=g:x\ \circ\vartriangle_{X\cdot Y}$.

Suppose we already proved $(i)^\circ$ and $(ii)^\circ$ for all ordinals $\iota<\lambda$ and $\varkappa<\mu$, where $\lambda$ and $\mu$ are two given ordinals:
\bit{
\item[$(i)^\circ$] the spaces $Y_\varkappa:X^\iota$ and $Y:X$ become isomorphic after pseudosaturation:
\beq\label{(Y_varkappa:X_iota)^vartriangle=(Y:X)^vartriangle}
(Y_\varkappa:X^\iota)^\vartriangle=(Y:X)^\vartriangle,
\eeq

\item[$(ii)^\circ$] for each $x\in X^\iota$ and $g\in (Y_\varkappa)^\star$ the functional
$$
    (g:x):(Y_\varkappa:X^\iota)\to\C\quad\Big|\quad (g:x)(\ph)=g(\ph(x)),\qquad \ph\in Y_\varkappa:X^\iota
$$
    is continuous on the space $(Y:X)^\vartriangle$, i.e. there exists a unique functional $h\in ((Y:X)^\vartriangle)^\star$ such that the following diagram is commutative:
\beq\label{(Y_varkappa:X_iota)^vartriangle=(Y:X)^vartriangle-func}
\xymatrix @R=2.pc @C=4.0pc 
{
(Y_\varkappa:X^\iota)^\vartriangle \ar[d]_{\vartriangle_{Y_\varkappa:X^\iota}} \ar[r]^{(\wedge_\varkappa^Y:\vee_\iota^X)^\vartriangle} & (Y:X)^\vartriangle \ar@{-->}[d]^{h} \\
Y_\varkappa:X^\iota \ar[r]_{g:x}& \C
}
\eeq
}\eit

Let us show then that (i) and (ii) hold for the ordinals $\lambda$ and $\mu$.

1. First we fix an arbitrary ordinal $\varkappa$ such that $\varkappa<\mu$, and show that $(i)^\circ$ and $(ii)^\circ$ hold after substitution $\iota=\lambda$. We have to consider two cases here.

a) Suppose $\lambda$ is an isolated ordinal, i.e.
$$
\lambda=\iota+1
$$
for some ordinal $\iota<\lambda$. Let us show that then $(ii)^\circ$ holds for $\lambda$ (substituted instead of $\iota$). Take
$$
x\in X^\lambda=X^{\iota+1}=(X^\iota)^\vee
$$
and
$$
g\in (Y_\varkappa)^\star
$$
Then $x$ is a limit of some totally bounded net $\{x_i;\ i\in I\}\subseteq X^\iota$:
\beq\label{(Y_mu:X_lambda)^vartriangle=(Y:X)^vartriangle-1}
x_i\overset{(X^\iota)^\vee}{\underset{i\to\infty}{\longrightarrow}}x
\eeq
Since the set $\{x_i\}\subseteq X^\iota$ is totally bounded, by Corollary \ref{COR:(g,x)|->g:x-razdelno-nepreryvno}, the set
$$
\{g:x_i;\ i\in I\}
$$
is totally bounded in $(Y_\varkappa:X^\iota)^\star$. On the other hand, from Theorem  \ref{TH:(g,x)|->g:x-razdelno-nepreryvno} we have that $\{g:x_i;\ i\in I\}$ is a Cauchy net. Thus, $\{g:x_i;\ i\in I\}$ is a totally bounded Cauchy net in $(Y_\varkappa:X^\iota)^\star$, and as a corollary, it converges in the enveloping pseudocomplete space $(Y_\varkappa:X^\iota)^{\star\triangledown}$:
$$
(Y_\varkappa:X^\iota)^\star\subseteq (Y_\varkappa:X^\iota)^{\star\triangledown}=
\text{\cite[Theorem 3.14]{Ak03}}= (Y_\varkappa:X^\iota)^{\vartriangle\star}= \eqref{(Y_varkappa:X_iota)^vartriangle=(Y:X)^vartriangle}=
(Y:X)^{\vartriangle\star}
$$
In other words, there is a functional $h\in (Y:X)^{\vartriangle\star}$ such that
\beq\label{(Y_mu:X_lambda)^vartriangle=(Y:X)^vartriangle-2}
g:x_i\overset{(Y:X)^{\vartriangle\star}}{\underset{i\to\infty}{\longrightarrow}}h
\eeq
On the other hand, from the chain of equalities
$$
(X^\iota)^\vee=X^{\iota+1}=X^\lambda
$$
we obtain that \eqref{(Y_mu:X_lambda)^vartriangle=(Y:X)^vartriangle-1} is equivalent to the relation
\beq\label{(Y_mu:X_lambda)^vartriangle=(Y:X)^vartriangle-3}
x_i\overset{X^\lambda}{\underset{i\to\infty}{\longrightarrow}} x
\eeq
which by Theorem \ref{TH:(g,x)|->g:x-razdelno-nepreryvno} implies
\beq\label{(Y_mu:X_lambda)^vartriangle=(Y:X)^vartriangle-4}
g:x_i\overset{(Y_\varkappa:X^\lambda)^{\star}}{\underset{i\to\infty}{\longrightarrow}} g:x
\eeq
Together \eqref{(Y_mu:X_lambda)^vartriangle=(Y:X)^vartriangle-2} and \eqref{(Y_mu:X_lambda)^vartriangle=(Y:X)^vartriangle-4} mean, in particular, that on each operator
$$
\ph\in Y_\varkappa:X^\lambda
$$
(by Lemma \ref{LM:Y_mu:X^lambda=Y:X}, the spaces $Y_\varkappa:X^\lambda$ and $Y:X$ coincide as sets, and the difference can be only in topologies) the functionals $g:x$ and $h$ coincide:
$$
(g:x)(\ph)=h(\ph).
$$
This proves $(ii)^\circ$ (for the case of $\lambda=\iota+1$ substituted instead of $\iota$).

Further we prove $(i)^\circ$ (again for the case of $\lambda=\iota+1$ substituted instead of $\iota$). The map
$$
(Y_\varkappa:X^\lambda)^\vartriangle\to (Y:X)^\vartriangle
$$
is continuous. We need to verify that its inverse map
$$
(Y:X)^\vartriangle\to (Y_\varkappa:X^\lambda)^\vartriangle.
$$
is continuous. Since the space $(Y:X)^\vartriangle$ is pseudosaturated, it is sufficient to prove the continuity of the map
$$
(Y:X)^\vartriangle\to (Y_\varkappa:X^\lambda).
$$
Take a basic neighbourhood of zero in $Y_\varkappa:X^\lambda$, i.e. a set $V:S$, where $V$ is a closed convex balanced neighbourhood of zero in $Y_\varkappa$, and $S$ a totally bounded set in $X^\lambda$. By Lemma \ref{LM:o-drobyah-v-Y_mu:X_lambda}, the set $V:S$ is capacious in $Y:X$. Hence if we prove that it is closed in  $(Y:X)^\vartriangle$, this will mean that it is a neighbourhood of zero in $(Y:X)^\vartriangle$.

This becomes obvious if we represent $V:S$ as the polar of the system of functionals of the form $\{g:x;\ g\in V^\circ,\ x\in S\}$:
\beq\label{(Y_mu:X_lambda)^vartriangle=(Y:X)^vartriangle-5}
V:S=\{g:x;\ g\in V^\circ,\ x\in S\}^\circ
\eeq
Indeed,
\begin{multline}\label{(Y_mu:X_lambda)^vartriangle=(Y:X)^vartriangle-5-proof}
\ph\in V:S \quad\Leftrightarrow\quad \ph(S)\subseteq V
\quad\Leftrightarrow\quad \forall x\in S \ \ph(x)\in V
\quad\Leftrightarrow \\ \Leftrightarrow\quad
\forall g\in V^\circ\ \forall x\in S \ |\underbrace{g(\ph(x))}_{\scriptsize \begin{matrix}\|\\ (g:x)(\ph)\end{matrix}}|\le 1
\quad\Leftrightarrow\quad \ph\in \{g:x;\ g\in V^\circ,\ x\in S\}^\circ
\end{multline}
Now we must notice that in formula \eqref{(Y_mu:X_lambda)^vartriangle=(Y:X)^vartriangle-5}
$x\in S$, where $S$ is a totally bounded set in $X^\lambda$, and $g\in V^\circ$, where $V$ is a neighbourhood of zero in $Y_\varkappa$. Hence, $x\in X^\lambda$, and $g\in (Y_\varkappa)^\star$. From this, by $(ii)^\circ$ (which we already proved for the case of $\lambda=\iota+1$), we have that $g:x$ are continuous functionals on the space $(Y:X)^\vartriangle$. We obtain, that $V:S$ is a polar of some system of continuous functionals on $(Y:X)^\vartriangle$, therefore, $V:S$ is a closed set in $(Y:X)^\vartriangle$. In addition it is capacious, hence it is a neighbourhood of zero in the (pseudosaturated) space $(Y:X)^\vartriangle$. We peroved $(i)^\circ$ (for the case of $\lambda=\iota+1$ substituted instead of $\iota$).

b) Now we consider the case when $\lambda$ is a limit ordinal, i.e.
$$
\lambda\ne \iota+1
$$
for all $\iota<\lambda$. Here we again must start with $(ii)^\circ$. Take $x\in X^\lambda$ and $g\in (Y_\varkappa)^\star$. Then
$$
x\in X^\lambda=\text{\cite[(1.19)]{Ak03}}=\bigcup_{\iota<\lambda} X^\iota,
$$
hence there exist an ordinal $\iota<\lambda$ such that $x\in X^\iota$. We obtain that $x\in X^\iota$ and $g\in (Y_\varkappa)^\star$, hence by assumption of the induction $(ii)^{\circ}$, the functional $g:x$ must be continuous on the space $(Y:X)^\vartriangle$. This proves $(ii)^\circ$ (for the case of the limit ordinal $\lambda$ substituted instead of $\iota<\lambda$). Now we can prove $(i)^\circ$. Here we repeat the reasoning of a): we have to prove the continuity of the mapping
$$
(Y:X)^\vartriangle\to (Y_\varkappa:X^\lambda).
$$
We take a basic neighbourhood of zero in $Y_\varkappa:X^\lambda$, i.e. a set $V:S$, where $S$ is totally bounded in $X^\lambda$, and $V$ a closed convex balanced neigbourood of zero in $Y_\varkappa$, then we prove equality
\eqref{(Y_mu:X_lambda)^vartriangle=(Y:X)^vartriangle-5}, then we notice that in this equality $x\in X^\lambda$, and  $g\in (Y_\varkappa)^\star$. After that by already proven $(ii)^\circ$ (for the case of the limit ordinal $\lambda$ substituted instead of $\iota<\lambda$), we obtain that $g:x$ are continuous functionals on the space $(Y:X)^\vartriangle$. Therefore, $V:S$ is a closed set in $(Y:X)^\vartriangle$. In addition it is capacious, by Lemma \ref{LM:o-drobyah-v-Y_mu:X_lambda}, hence it is a neighbourhood of zero in the  (pseudosaturated) space $(Y:X)^\vartriangle$. We proved $(i)^\circ$ (for the case of the limit ordinal $\lambda$ is substituted instead of $\iota<\lambda$).

2. We showed that for a given $\varkappa<\mu$ the conditions $(i)^\circ$ and $(ii)^\circ$ hold if we replace there $\iota<\lambda$ by $\iota=\lambda$. Now let us take $\iota=\lambda$, and show that $(i)^\circ$ and $(ii)^\circ$ (with  $\iota=\lambda$) remain true if we replace there $\varkappa<\mu$ by $\varkappa=\mu$. Again we have to consider two cases.

a) We assume first that $\mu$ is an isolated ordinal, i.e.
$$
\mu=\varkappa+1
$$
for some ordinal $\varkappa<\mu$. Let us show that for $\mu$ (substituted instead of $\varkappa$) the condition $(ii)^\circ$ holds as well. Take
$$
g\in (Y_\mu)^\star=(Y_{\varkappa+1})^\star=((Y_\varkappa)^\wedge)^\star
=\text{\cite[Theorem 3.10]{Ak03}}=((Y_\varkappa)^\star)^\vee
$$
and
$$
x\in X^\lambda.
$$
Then $g$ is a limit of some totally bounded net $\{g_i;\ i\in I\}\subseteq (Y_\varkappa)^\star$:
\beq\label{(Y_mu:X_lambda)^vartriangle=(Y:X)^vartriangle-1*}
g_i\overset{((Y_\varkappa)^\star)^\vee}{\underset{i\to\infty}{\longrightarrow}}g
\eeq
Since the set $\{g_i\}\subseteq (Y_\varkappa)^\star$ is totally bounded, by Corollary \ref{COR:(g,x)|->g:x-razdelno-nepreryvno}, the set
$$
\{g_i:x;\ i\in I\}
$$
is totally bounded in $(Y_\varkappa:X^\lambda)^\star$. On the other hand, Theorem \ref{TH:(g,x)|->g:x-razdelno-nepreryvno} implies that $\{g_i:x;\ i\in I\}$ is a Cauchy net. Thus, we obtain that $\{g_i:x;\ i\in I\}$ is a totally bounded Cauchy net in the space $(Y_\varkappa:X^\lambda)^\star$, hence, it converges in the enveloping pseudocomplete space $(Y_\varkappa:X^\lambda)^{\star\triangledown}$:
$$
(Y_\varkappa:X^\lambda)^\star\subseteq (Y_\varkappa:X^\lambda)^{\star\triangledown}=
\text{\cite[Theorem 3.14]{Ak03}}= (Y_\varkappa:X^\lambda)^{\vartriangle\star}= \eqref{(Y_varkappa:X_iota)^vartriangle=(Y:X)^vartriangle}=
(Y:X)^{\vartriangle\star}
$$
I.e. there exists a functional $h\in (Y:X)^{\vartriangle\star}$ such that
\beq\label{(Y_mu:X_lambda)^vartriangle=(Y:X)^vartriangle-2*}
g_i:x\overset{(Y:X)^{\vartriangle\star}}{\underset{i\to\infty}{\longrightarrow}}h
\eeq
On the other hand, from the chain of equalities
$$
((Y_\varkappa)^\star)^\vee= \text{\cite[Theorem 3.10]{Ak03}}= ((Y_\varkappa)^\wedge)^\star=(Y_{\varkappa+1})^\star=(Y_\mu)^\star
$$
we have that the relation \eqref{(Y_mu:X_lambda)^vartriangle=(Y:X)^vartriangle-1*} is equivalent to the relation
\beq\label{(Y_mu:X_lambda)^vartriangle=(Y:X)^vartriangle-3*}
g_i\overset{(Y_\mu)^\star}{\underset{i\to\infty}{\longrightarrow}} g
\eeq
from which by Theorem \ref{TH:(g,x)|->g:x-razdelno-nepreryvno} it follows that
\beq\label{(Y_mu:X_lambda)^vartriangle=(Y:X)^vartriangle-4*}
g_i:x\overset{(Y_\mu:X^\lambda)^{\star}}{\underset{i\to\infty}{\longrightarrow}} g:x
\eeq
Together \eqref{(Y_mu:X_lambda)^vartriangle=(Y:X)^vartriangle-2*} and \eqref{(Y_mu:X_lambda)^vartriangle=(Y:X)^vartriangle-4*} mean, in particular, that on each operator
$$
\ph\in Y_\mu:X^\lambda
$$
(by Lemma \ref{LM:Y_mu:X^lambda=Y:X}, the spaces $Y_\mu:X^\lambda$ and $Y:X$ coincide as sets, but the difference can be in topologies) the functionals $g:x$ and $h$ coincide:
$$
(g:x)(\ph)=h(\ph).
$$
This proves $(ii)^\circ$ (for the case $\mu=\varkappa+1$ substituted instead of $\varkappa$).

Let us prove now $(i)^\circ$ (again for the case $\mu=\varkappa+1$ substituted instead of $\varkappa$). The mapping
$$
(Y_\mu:X^\lambda)^\vartriangle\to (Y:X)^\vartriangle
$$
is continuous. We need to prove that its inverse mapping
$$
(Y:X)^\vartriangle\to (Y_\mu:X^\lambda)^\vartriangle
$$
is continuous. Since the space $(Y:X)^\vartriangle$ is pseudosaturated, it is sufficient to prove the continuity of the mapping
$$
(Y:X)^\vartriangle\to (Y_\mu:X^\lambda).
$$
Consider the basic neighbourhood of zero in $Y_\mu:X^\lambda$, i.e. a set $V:S$, where $V$ is a closed convex balanced nieghbourhood of zero in $Y_\mu$, and $S$ a totally bounded set in $X^\lambda$. By Lemma \ref{LM:o-drobyah-v-Y_mu:X_lambda}, the set $V:S$ is capacious in $Y:X$. Hence if we prove that it is closed in $(Y:X)^\vartriangle$, this will mean that it is a neighbourhood of zero in $(Y:X)^\vartriangle$.

This becomes obvious if we represent $V:S$ as a polar of the system of functionals of the form $\{g:x;\ g\in V^\circ,\ x\in S\}$:
\beq\label{(Y_mu:X_lambda)^vartriangle=(Y:X)^vartriangle-5*}
V:S=\{g:x;\ g\in V^\circ,\ x\in S\}^\circ
\eeq
which is proved by the same chain \eqref{(Y_mu:X_lambda)^vartriangle=(Y:X)^vartriangle-5-proof} as before.
Further we notice that in formula \eqref{(Y_mu:X_lambda)^vartriangle=(Y:X)^vartriangle-5*}
$x\in S$, where $S$ is a totally bounded set in $X^\lambda$, and $g\in V^\circ$, where $V$ is a neighbourhood of zero in $Y_\mu$. Hence, $x\in X^\lambda$, and $g\in (Y_\mu)^\star$. From this, by $(ii)^\circ$ (which we already proved for the case $\mu=\varkappa+1$),  $g:x$ are continuous functionals on the space $(Y:X)^\vartriangle$. We see that $V:S$ is a polar of a system of continuous functionals on $(Y:X)^\vartriangle$, hence $V:S$ is a closed set in  $(Y:X)^\vartriangle$. In addition, it is capacious, hence it is a neighbourhood of zero in the (pseudosaturated) space $(Y:X)^\vartriangle$. We proved $(i)^\circ$ (for the case $\mu=\varkappa+1$ substituted instead of $\varkappa$).

b) Now let us consider the case where $\mu$ is a limit ordinal, i.e.
$$
\mu\ne \varkappa+1
$$
for all $\varkappa<\mu$. Here we again start with $(ii)^\circ$. Take $x\in X^\lambda$ and $g\in (Y_\mu)^\star$. Then
$$
g\in (Y_\mu)^\star=\text{\cite[Theorem 3.12]{Ak03}} =(Y^\star)^\mu=\text{\cite[(1.19)]{Ak03}}=\bigcup_{\varkappa<\mu} (Y^\star)^\varkappa=\text{\cite[Theorem 3.12]{Ak03}} =\bigcup_{\varkappa<\mu}(Y_\varkappa)^\star,
$$
so there is an ordinal $\varkappa<\mu$ such that $g\in (Y^\star)^\varkappa=(Y_\varkappa)^\star$. We have that $x\in X^\iota$ and $g\in (Y_\varkappa)^\star$, hence by the assumption of induction in $(ii)^{\circ}$, the functional $g:x$ must be continuous on the space $(Y:X)^\vartriangle$. This proves  $(ii)^\circ$ (for the case of limit $\mu$ substituted instead of $\varkappa$). We can now prove $(i)^\circ$. We repeat here the reasoning of a): we need to prove the continuity of the mapping
$$
(Y:X)^\vartriangle\to (Y_\mu:X^\lambda).
$$
We take a basic neighbourhood of zero in $Y_\mu:X^\lambda$, i.e. a set $V:S$, where  $S$ is a totally bounded set in $X^\lambda$, and $V$ a closed convex balanced neighbourhood of zero in $Y_\mu$, then we prove
\eqref{(Y_mu:X_lambda)^vartriangle=(Y:X)^vartriangle-5*}, then we notice that there $x\in X^\lambda$, and $g\in (Y_\mu)^\star$. After that by $(ii)^\circ$ (which we already proved for limit $\mu$ substituted instead of $\varkappa$), we obtain that $g:x$ are continuous functional on the space $(Y:X)^\vartriangle$. Hence $V:S$ is a closed set in $(Y:X)^\vartriangle$. In addition it is capacious by Lemma \ref{LM:o-drobyah-v-Y_mu:X_lambda}, therefore it is a neighbourhood of zero in the (pseudosaturated) space $(Y:X)^\vartriangle$. We proved $(i)^\circ$ (for the case of limit $\mu$ substituted instead of $\varkappa$).
\epr

For arbitrary locally convex spaces $X$ and $Y$ the mappings
$$
\triangledown_X: X\to X^\triangledown,\qquad \vartriangle_Y: Y^\vartriangle\to Y
$$
define the mapping
\beq\label{(vartriangle_Y:triangledown_X)^vartriangle}
(\vartriangle_Y:\triangledown_X): (Y^\vartriangle:X^\triangledown) \to (Y:X).
\eeq
(which, certainly, coincides with the mapping \eqref{wedge_mu^Y:vee_lambda^X-1} for some $\lambda$ and $\mu$).

\btm\label{TH:(vartriangle_Y:triangledown_X)^vartriangle-isomorphism}
Let $X$ be pseudosaturated and $Y$ pseudocomplete. Then the pseudosaturation of the mapping \eqref{(vartriangle_Y:triangledown_X)^vartriangle} is an isomorphism of locally convex (and stereotype) spaces:
\beq\label{(vartriangle_Y:triangledown_X)^vartriangle-isomorphism}
Y^\vartriangle\oslash X^\triangledown\cong (Y^\vartriangle:X^\triangledown)^\vartriangle\cong  (Y:X)^\vartriangle
\eeq
\etm
\bpr
The second equality follows from \eqref{(Y_mu:X_lambda)^vartriangle=(Y:X)^vartriangle}, if we choose ordinals $\lambda$ and $\mu$ such that $X^\lambda=X^\triangledown$ and $Y_\mu=Y^\vartriangle$.
And the first one is just definition \eqref{DEF:oslash} of the tensor fraction $\oslash$ (here the spaces $X^\triangledown$ and $Y^\vartriangle$ are stereotype by \cite[Proposition 3.17]{Ak03} and \cite[Proposition 3.16]{Ak03}).
\epr

\section{Pseudocompletion of the primary tensor product $X\cdot Y$}

\subsection{Primary tensor product $X\cdot Y$.}

Let $X$ and $Y$ be locally convex spaces. Let us call their {\it primary tensor product} $X\cdot Y$ the locally convex space, consisting of linear continuous mappings $\ph:X^\star\to Y$, endowed with the {\it topology of uniform convergence on polars of the neighbourhoods of zero $U\subseteq X$}:
\beq\label{shodimost-v-X-cdot-Y}
\ph_i\overset{X\cdot Y}{\underset{i\to\infty}{\longrightarrow}}\ph\quad\Leftrightarrow\quad
\forall U\in{\mathcal U}(X)\quad
\ph_i(f)\overset{Y}{\underset{f\in U^\circ}{\underset{i\to\infty}{\rightrightarrows}}}\ph(f)
\eeq
It will be useful to choose a letter for this topology, let it be $\xi$. Then the space $X\cdot Y$ can be presented by the formula
\beq\label{DEF:X-cdot-Y}
X\cdot Y=Y\underset{\xi}{:} X^\star
\eeq
(the index $\xi$ means the convergence in the topology $\xi$). Obviously, there is a bijective linear continuous mapping
$$
Y:X\to Y\underset{\xi}{:} X^\star=X\cdot Y
$$
(which however is not an isomorphism if $X$ is not pseudosaturated).
If $A\subseteq X$ and $B\subseteq Y$, then the symbol $A\cdot B$ denotes the set
$$
A\cdot B=B:(A^\circ)
$$
(where the colon : was defined in \eqref{DEF:B:A}).

If $\alpha:X'\to X$, $\beta:Y'\to Y$ are linear continuous mappings, then the diagram
\beq\label{alpha-cdot-beta}
\xymatrix @R=2.pc @C=2.0pc 
{
X^\star \ar[d]_{(\alpha\cdot\beta)(\psi)=\beta\circ\psi\circ\alpha^\star} \ar[r]^{\alpha^\star} & (X')^\star \ar[d]^{\psi} \\
Y & Y' \ar[l]_{\beta}
}
\eeq
defines a mapping
$$
\ph\cdot\chi: X'\cdot Y'\to X\cdot Y.
$$

\bprop\label{alpha-cdot-beta-nepreryvnoe-otobr}
For any operators $\alpha:X'\to X$ and $\beta:Y'\to Y$ the mapping
$$
\ph\cdot\chi: X'\cdot Y'\to X\cdot Y
$$
is continuous.
\eprop
\bpr The continuity of $\alpha$ implies that the dual operator  $\alpha^\star:X^\star\to (X')^\star$ turns the polar of arbitrary neighbourhood of zero $U\subseteq X$ into the polar of some neighbourhood of zero $U'$, namely, $U'=\alpha^{-1}(U)$:
\begin{multline*}
f\in U^\circ\quad\Rightarrow\quad \forall x\in U\ \abs{f(x)}\le 1
\quad\Rightarrow\\ \Rightarrow\quad \forall x'\in U'=\alpha^{-1}(U)\ \abs{\alpha^\star(f)(x')}=|f(\underbrace{\alpha(x')}_{\scriptsize \begin{matrix}\text{\rotatebox{90}{$\owns$}}\\ U\end{matrix}})|\le 1
\quad\Rightarrow\quad \alpha^\star(f)\in (U')^\circ
\end{multline*}
As a corollary the convergence $\psi_i\to \psi$ in $X'\cdot Y'$
$$
\forall U'\in{\mathcal U}(X')\quad
\psi_i(g)\overset{Y'}{\underset{g\in (U')^\circ}{\underset{i\to\infty}{\rightrightarrows}}}\psi(g)
$$
implies the convergence $(\alpha\cdot\beta)(\psi_i)\to (\alpha\cdot\beta)(\psi)$ in $X\cdot Y$:
$$
\forall U\in{\mathcal U}(X)\quad
(\alpha\cdot\beta)(\psi_i)(f)=\beta\Big(\psi_i\big(\alpha^\star(f)\big)\Big)
\overset{Y}{\underset{f\in U^\circ}{\underset{i\to\infty}{\rightrightarrows}}}
\beta\Big(\psi\big(\alpha^\star(f)\big)\Big)=
(\alpha\cdot\beta)(\psi)(f)
$$
\epr

Suppose further that $X$ and $Y$ are locally convex spaces, and $X$ is pseudocomplete. For each operator $\ph:X^\star\to Y$ we can consider the dual operator $\ph:Y^\star\to X^{\star\star}$, and after that form the composition with the operator $i_X^{-1}:X^{\star\star}\to X$ (which due to \cite[Corollary 2.13]{Ak03} exists and is continuous, since $X$ is pseudocomplete). Let us denote this composition by $\omega_{X,Y}(\ph)$:
\beq\label{DEF:omega(ph)}
\xymatrix @R=2.pc @C=6.0pc 
{
Y^\star \ar@/_4ex/[rr]_{\omega_{X,Y}(\ph)} \ar[r]^{\ph^\star} & X^{\star\star} \ar[r]^{i_X^{-1}} & X
}
\eeq
We have the mapping
\beq\label{omega:XY->YX}
\omega_{X,Y}:X\cdot Y=Y\underset{\xi}{:}X^\star\to X\underset{\xi}{:}Y^\star=Y\cdot X\quad\Big|\quad \omega_{X,Y}(\ph)=i_X^{-1}\circ\ph^\star
\eeq

\btm\label{TH:X-cdot-Y-cong-Y-cdot-X}
If locally convex spaces $X$ and $Y$ are pseudocomplete, then the mapping \eqref{omega:XY->YX} establishes an isomorphism of locally convex spaces
\beq\label{X-cdot-Y-cong-Y-cdot-X}
X\cdot Y\cong Y\cdot X
\eeq
\etm
\bpr
1. First we must notice the following identity:
\beq\label{f(i_X^-1(ph^star(g)))=g((ph(f)))}
f(\omega_{X,Y}(\ph)(g))=g(\ph(f)),\qquad f\in X^\star, \ g\in Y^\star,\ \ph\in Y:(X^\star).
\eeq
For proof we put
$$
x=\omega_{X,Y}(\ph)(g)=i_X^{-1}(\ph^\star(g))=i_X^{-1}(g\circ\ph)\in X,
$$
then we have
$$
f(\underbrace{\omega_{X,Y}(\ph)(g)}_{\scriptsize\begin{matrix}\|\\ x\end{matrix}})=f(x)=i_X(x)(f)=
i_X(i_X^{-1}(g\circ\ph))(f)=(g\circ\ph)(f)=g(\ph(f))
$$

2. When \eqref{f(i_X^-1(ph^star(g)))=g((ph(f)))} is proved, we must verify that the mapping $\omega_{X,Y}$ is bijective. This follows from the identity
\beq
\omega_{Y,X}(\omega_{X,Y}(\ph))=\ph,\qquad \ph\in X\cdot Y=Y\underset{\xi}{:}X^\star.
\eeq
Indeed, for each $f\in X^\star$ and $g\in Y^\star$ we have
$$
g(\omega_{Y,X}(\omega_{X,Y}(\ph))(f))=\eqref{f(i_X^-1(ph^star(g)))=g((ph(f)))}=
f(\omega_{X,Y}(\ph)(g))=\eqref{f(i_X^-1(ph^star(g)))=g((ph(f)))}=g(\ph(f)).
$$

3. Finally, let us prove the continuity of the mapping $\omega_{X,Y}$ in both directions. Indeed, if $U:(V^\circ)$ is a basic neighbourhood of zero in $X\underset{\xi}{:}Y^\star=Y\cdot X$ (where $U\subseteq X$ and $V\subseteq Y$ are closed convex balanced neighbourhood of zero), then $V:(U^\circ)$ is a basic neighbourhood of zero in $Y\underset{\xi}{:}X^\star=X\cdot Y$, and they turn into each other under the mappings $\omega_{X,Y}$ and $\omega_{X,Y}^{-1}$:
\begin{multline*}
\omega_{X,Y}(\ph)\in U:(V^\circ)\quad\Leftrightarrow\quad \omega_{X,Y}(\ph)(V^\circ)\subseteq U
\quad\Leftrightarrow\quad \forall f\in U^\circ\ \forall g\in V^\circ\quad |\underbrace{f(\omega_{X,Y}(\ph)(g))}_{\scriptsize \begin{matrix}\|\put(3,0){\eqref{f(i_X^-1(ph^star(g)))=g((ph(f)))}}
\\ g(\ph(f))\end{matrix}}|\le 1
\quad\Leftrightarrow
\\ \Leftrightarrow \quad
\forall f\in U^\circ\ \forall g\in V^\circ\quad |g(\ph(f))|\le 1
\quad\Leftrightarrow\quad \ph(U^\circ)\subseteq V \quad\Leftrightarrow\quad \ph\in V:(U^\circ)
\end{multline*}
\epr

Let us go back to diagram \eqref{alpha-cdot-beta}. In the special case when $\alpha$ and $\beta$ are functionals
$$
f:X\to\C,\qquad g:Y\to\C,
$$
the mapping $\alpha\cdot\beta$ can be naturally understood as the functional
$$
f\cdot g:X\cdot Y\to\C,
$$
acting by the formula
\beq\label{DEF:f-cdot-g}
(f\cdot g)(\ph)=g(\ph(f)),\qquad \ph:X^\star\to Y.
\eeq

If $F\subseteq X^\star$ and $G\subseteq Y^\star$ are two sets of functionals, then the symbol $F\cdot G$ denotes the set
$$
F\cdot G=\{f\cdot g;\ f\in F,\ g\in G\}.
$$
In particular, if $f\in X^\star$ and $G\subseteq Y^\star$, then
$$
f\cdot G=\{f\cdot g;\ g\in G\}.
$$
and if $F\subseteq X^\star$ and $g\in Y^\star$, then
$$
F\cdot g=\{f\cdot g;\ f\in F\}.
$$

\btm\label{TH:(f,g)|->f-cdot-g-razdelno-nepreryvno}
For any pseudocomplete locally convex spaces $X$ and $Y$ the mapping
$$
(f,g)\in X^\star\times Y^\star\mapsto f\cdot g\in (X\cdot Y)^\star
$$
is separately continuous.
\etm
\bpr
Let $\varPhi\subseteq X\cdot Y=Y\underset{\xi}{:}X^\star$ be a totally bounded set, and $f\in X^\star$. By definition of the topology $\xi$, $\varPhi$ is equicontinuous and uniformly totally bounded on polar $U^\circ$ of arbitrary neighbourhood of zero $U\subseteq X$ (we use here a variant of \cite[Theorem 5.1]{Ak03}). In particular, if we take the neighbourhood $U=\{f\}^\circ$, then the set of operators $\varPhi$ turns its polar
$$
U^\circ=\{f\}^{\circ\circ}=\{\lambda\in\C:\ |\lambda|\le 1\}\cdot f.
$$
into a totally bounded set in $Y$:
$$
\varPhi(U^\circ)=\{\lambda\in\C:\ |\lambda|\le 1\}\cdot \varPhi(f)\subseteq Y.
$$
As a corollary, the polar of this set
$$
V=\{\lambda\in\C:\ |\lambda|\le 1\}\cdot \varPhi(f)^\circ\subseteq Y^\star
$$
is a neighbourhood of zero in $Y^\star$. For each $g\in V$ we have
$$
\sup_{\ph\in\varPhi}\abs{(f\cdot g)(\ph)}=\sup_{\ph\in\varPhi}\abs{g(\ph(f))}\le 1
$$
i.e. $f\cdot g\in\varPhi^\circ$. In other words,
\beq\label{(f,g)|->f-cdot-g-razdelno-nepreryvno-1}
f\cdot V\subseteq \varPhi^\circ.
\eeq

We can say that for each basic neighbourhood of zero $\varPhi^\circ$ in $(X\cdot Y)^\star$ (where $\varPhi$ is a totally bounded set in $X\cdot Y$) and for any vector $f\in X^\star$ there is a neighbourhood of zero $V\subseteq Y^\star$ such that \eqref{(f,g)|->f-cdot-g-razdelno-nepreryvno-1} holds. By \eqref{X-cdot-Y-cong-Y-cdot-X}, the same is true if we transpose $X$ and $Y$.
\epr

\bcor\label{LM:F-cdot-g-subseteq-(X-cdot-Y)^star}
If $F\subseteq X^\star$ is totally bounded, and $g\in Y^\star$, then $F\cdot g\subseteq (X\cdot Y)^\star$ is totally bounded. Dually, if $f\in X^\star$ and $G\subseteq Y^\star$ is totally bounded, then $f\cdot G\subseteq (X\cdot Y)^\star$ is totally bounded.
\ecor

\bex\label{PROP:X-cdot-Y-cong-X-widetilde-otimes_e-Y}
Let $X$ and $Y$ be complete locally convex spaces, with $Y$ having the (classical) approximation property. Then the primary tensor product $X\cdot Y$ is isomorphic to the injective tensor producst
\beq\label{X-cdot-Y-cong-X-widetilde-otimes_e-Y}
X\cdot Y\cong X\widetilde{\otimes}_{\e} Y.
\eeq
If in addition $Y$ is nuclear, then $X\cdot Y$ is isomorphic to the projective tensor product
\beq\label{X-cdot-Y-cong-X-widetilde-otimes_pi-Y}
X\cdot Y\cong X\widetilde{\otimes}_{\e} Y\cong X\widetilde{\otimes}_{\pi} Y.
\eeq
\eex
\bpr
1. Consider the embedding $X\otimes Y\subseteq X\cdot Y=Y\underset{\xi}{:} X^\star$, defined by the formula
$$
(x\otimes y)(f)=f(x)\cdot y,\qquad x\in X,\ y\in Y, \ f\in X^\star.
$$
The topology $\xi$ on $X\otimes Y$, induced from $Y\underset{\xi}{:} X^\star$, coincides with the injective topology $\e$, since $\xi$ and $\e$ are generated by one system of seminorms:
$$
\abs{z}_{U,V}=\sup_{f\in U^\circ,\ g\in V^\circ}\abs{g(z(f))}=
\sup_{f\in U^\circ,\ g\in V^\circ}\abs{g\bigg(\Big(\sum_{i=1}^n x_i\otimes y_i\Big)(f)\bigg)}=
\sup_{f\in U^\circ,\ g\in V^\circ}\abs{\sum_{i=1}^n f(x_i)\cdot g(y_i)}
$$
for $z=\sum_{i=1}^n x_i\otimes y_i\in X\otimes Y$.

On the other hand, the space $Y\underset{\xi}{:} X^\star$ is complete, since $Y$ is complete, and $X^\star$ is saturated (as the dual to the complete space $X$, by \cite[Theorem 2.17]{Ak03}). This implies that the injective tensor product $X\widetilde{\otimes}_{\e} Y$ (i.e. the completion of $X\otimes Y$ with respect to the topology $\e$) is just the closure of the space $X\otimes Y$ in the space $X\cdot Y=Y\underset{\xi}{:} X^\star$:
$$
X\widetilde{\otimes}_{\e} Y\cong\overline{X\otimes Y}.
$$
To prove \eqref{X-cdot-Y-cong-X-widetilde-otimes_e-Y} we need only to verufy that $X\otimes Y$ is dense in $X\cdot Y=Y\underset{\xi}{:} X^\star$:
$$
\overline{X\otimes Y}=X\cdot Y.
$$
This follows from the fact that $Y$ has the approximation property: if $\psi_i$ is a net of operators of finite rank, tending to the identity operator uniformly on compact sets in $Y$,
$$
\psi_i\underset{i\to\infty}{\longrightarrow}\id_Y
$$
then for each $\ph\in X\cdot Y=Y\underset{\xi}{:} X^\star$ we have
$$
\psi_i\circ\ph\overset{X\cdot Y}{\underset{i\to\infty}{\longrightarrow}}\id_Y\circ\ph=\ph.
$$

2. If in addition $Y$ is nuclear, then the second identity in \eqref{X-cdot-Y-cong-X-widetilde-otimes_pi-Y} comes automatically from the characterization of nuclear spaces \cite[7.3.3]{Pietsch}:
$$
X\widetilde{\otimes}_{\e} Y\cong X\widetilde{\otimes}_{\pi} Y.
$$
\epr

\subsection{Pseudosaturation of $X\cdot Y$.}

Let $\{X_\lambda;\ \lambda\in\bf{Ord}\}$ and $\{Y_\mu;\ \mu\in\bf{Ord}\}$ be projective series of the spaces $X$ and $Y$ \cite[1.4(b)]{Ak03}, and
$$
\wedge_\lambda^X: X_\lambda\to X,\qquad \wedge_\mu^Y: Y_\mu\to Y
$$
the corresponding natural mappings. If we consider the analog of the diagram \eqref{alpha-cdot-beta}
\beq\label{wedge_lambda^X-cdot-wedge_mu^Y}
\xymatrix @R=2.pc @C=2.0pc 
{
X^\star \ar[d]_{(\wedge_\lambda^X\cdot\wedge_\mu^Y)(\psi)= \wedge_\mu^Y\circ\psi\circ(\wedge_\lambda^X)^\star} \ar[r]^{(\wedge_\lambda^X)^\star} & (X_\lambda)^\star \ar[d]^{\psi} \\
Y & Y_\mu \ar[l]_{\wedge_\mu^Y}
}
\eeq
we obtain the mapping
\beq\label{wedge_lambda^X-cdot-wedge_mu^Y-1}
\wedge_\lambda^X\cdot\wedge_\mu^Y: X_\lambda\cdot Y_\mu\to X\cdot Y,
\eeq
which is continuous by Proposition \ref{alpha-cdot-beta-nepreryvnoe-otobr}.

\blm\label{LM:o-drobyah-*}
If $X$ and $Y$ are pseudocomplete, then the mapping \eqref{wedge_lambda^X-cdot-wedge_mu^Y-1} is a bijection, hence the spaces $X_\lambda\cdot Y_\mu$ and $X\cdot Y$ can be identified as sets (with the possible difference in topologies):
\beq\label{X_lambda-cdot-Y_mu=X-cdot-Y}
X_\lambda\cdot Y_\mu=X\cdot Y
\eeq
\elm
\bpr
The mapping $\wedge_\lambda^X: X_\lambda\to X$ is an injection, hence that dual mapping $(\wedge_\lambda^X)^\star: (X_\lambda)^\star\gets X^\star$ is an epimorphism. On the other hand, the mapping $\wedge_\mu^Y: Y_\mu\to Y$ is an injection. This implies that the mapping
$$
\psi\mapsto \wedge_\mu^Y\circ\psi\circ (\wedge_\lambda^X)^\star
$$
is an injection. We only have to verify that it is a surjection. Indeed, let $\ph:X^\star\to Y$ be an arbitrary operator. We subsequently construct two operators $\ph_\mu$ and $\psi$, such that the following diagram is commutative:
$$
\xymatrix @R=3.pc @C=4.0pc 
{
X^\star \ar[d]_{\ph}\ar@{-->}[dr]_{\ph_\mu} \ar[r]^{(\wedge_\lambda^X)^\star} & (X_\lambda)^\star \ar@{-->}[d]^{\psi} \\
Y & Y_\mu \ar[l]^{\wedge_\mu^Y}
}
$$
Since $X$ is pseudocomplete, the space $X^\star$ is pseudosaturated. This means that $\ph$ is uniquely extended to an operator $\ph^\vartriangle:X^\star\to Y^\vartriangle$. If we consider its composition with the natural mapping $Y^\vartriangle\to Y_\mu$, we obtain an operator, extending $\ph$ to an operator $\ph_\mu:X^\star\to Y_\mu$.

Further, since $Y$ is pseudocomplete, $Y_\mu$ is also pseudocomplete. Hence the operator $\ph_\mu:X^\star\to Y_\mu$ is extended to some operator $\chi:(X^\star)^\triangledown\to Y_\mu$. If we consider its composition with the natural embedding $(X^\star)^\lambda\to (X^\star)^\triangledown$, we obtain an operator $\psi:(X^\star)^\lambda\to Y_\mu$. The equality \cite[(3.11)]{Ak03}
$$
(X^\star)^\lambda=(X_\lambda)^\star
$$
allows to present $\psi$ as an operator $\psi:(X_\lambda)^\star\to Y_\mu$.
\epr

\blm\label{LM:o-sohranenii-polnoi-ogranichennosti}
Let $X$ and $Y$ be pseudocomplete. If $\varPhi\in X\cdot Y$ is a totally bounded set, then for any $\lambda,\mu\in\bf{Ord}$ its representation in the space $X_\lambda\cdot Y_\mu$ (by the bijection \eqref{X_lambda-cdot-Y_mu=X-cdot-Y}) is also a totally bounded set.
\elm
\bpr
1. First, consider the case, when $\lambda=0$. If $\varPhi\subseteq Y\underset{\xi}{:}X^\star$ is a totally bounded set, then this means that $\varPhi$ is equicontinuous and uniformly totally bounded on the polar $U^\circ$ of every neighbourhood of zero $U\subseteq X$ (we use here a variant of \cite[Theorem 5.1]{Ak03}). As a corollary, the image $\varPhi(U^\circ)$ is a totally bounded subset in $Y$. The space $Y_\mu$ is a strengthening of the topology on $Y$, in which the class of totally bounded sets and the topology on these sets are not changed. Hence $\varPhi(U^\circ)$ is totally bounded in the space $Y_\mu$ as well, and moreover, with the same induced topology as the one induced from $Y$. We can conclude from this that $\varPhi$, considered as a set of operators from $X^\star$ to $Y_\mu$, is equicontinuous and uniformly totally bounded on $U^\circ$. And this is true for each neighbourhood of zero $U\subseteq X$. Hence, $\varPhi$ is totally bounded in $Y_\mu\underset{\xi}{:}X^\star=X\cdot Y_\mu$.

2. Thus we understood that if $\varPhi\subseteq X\cdot Y$ is totally bounded in $X\cdot Y$, then it is totally bounded in $X\cdot Y_\mu$ for each $\mu\in\bf{Ord}$. By Theorem \ref{TH:X-cdot-Y-cong-Y-cdot-X} this means that $\varPhi$ is totally bounded in the space $Y_\mu\cdot X\cong X\cdot Y_\mu$. Again what we have already proved allows us to conclude that $\varPhi$ is totally bounded in the space $Y_\mu\cdot X_\lambda$ for each $\lambda\in\bf{Ord}$. And again applying Theorem \ref{TH:X-cdot-Y-cong-Y-cdot-X}, we can say that $\varPhi$ is totally bounded in the space $X_\lambda\cdot Y_\mu\cong Y_\mu\cdot X_\lambda$.
\epr

\blm\label{LM:o-drobyah}
Let $X$ and $Y$ be pseudocomplete. If $U\in {\mathcal U}(X_\lambda)$ and $V\in {\mathcal U}(Y_\mu)$, then $U\cdot V\in {\mathcal D}(X\cdot Y)$.
\elm
\bpr
Let $\varPhi\subseteq X\cdot Y=Y\underset{\xi}{:}X^\star$ be a totally bounded set. By Lemma \ref{LM:o-sohranenii-polnoi-ogranichennosti}, $\varPhi$ is totally bounded in the space $X_\lambda\cdot Y_\mu=Y_\mu\underset{\xi}{:}(X_\lambda)^\star$. Since $U\cdot V=U:(V^\circ)$ is a neighbourhood of zero in the space $X_\lambda\cdot Y_\mu=Y_\mu\underset{\xi}{:}(X_\lambda)^\star$, there is a finite set $A\subseteq X_\lambda\cdot Y_\mu=Y_\mu\underset{\xi}{:}(X_\lambda)^\star$ such that
\beq\label{o-drobyah-1}
\varPhi\subseteq U\cdot V+A.
\eeq
But by Lemma \ref{LM:o-drobyah-*} we can think that all these sets are contained in the space $X\cdot Y$. We obtain the following: for any totally bounded set $\varPhi$ in $X\cdot Y$ there exists a finite set $A$ in $X\cdot Y$ such that \eqref{o-drobyah-1} holds. This means that the set $U\cdot V$ is capacious in the space $X\cdot Y$.
\epr

\blm\label{LM:(X_lambda-cdot-Y_lambda)^vartriangle=(X-cdot-Y)^vartriangle}
For any pseudocomplete locally convex space $X$ and $Y$ and for any ordinals $\lambda,\mu\in\bf{Ord}$
\bit{
\item[(i)] the spaces $X_\lambda\cdot Y_\mu$ and $X\cdot Y$ become isomorphic after pseudosaturation:
\beq\label{(X_lambda-cdot-Y_mu)^vartriangle=(X-cdot-Y)^vartriangle}
(X_\lambda\cdot Y_\mu)^\vartriangle=(X\cdot Y)^\vartriangle,
\eeq

\item[(ii)] for any functionals $f\in (X_\lambda)^\star$ and $g\in (Y_\mu)^\star$ the functional $f\cdot g:X_\lambda\cdot Y_\mu\to\C$ is continuous on the space $(X\cdot Y)^\vartriangle$, i.e. there exists a unique functional $h\in ((X\cdot Y)^\vartriangle)^\star$, such that the following diagram is commutative:
\beq\label{wedge_lambda^X-cdot-wedge_mu^Y}
\xymatrix @R=2.pc @C=4.0pc 
{
(X_\lambda\cdot Y_\mu)^\vartriangle \ar[d]_{\vartriangle_{X_\lambda\cdot Y_\mu}} \ar[r]^{(\wedge_\lambda^X\cdot\wedge_\mu^Y)^\vartriangle} & (X\cdot Y)^\vartriangle \ar@{-->}[d]^{h} \\
X_\lambda\cdot Y_\mu \ar[r]_{f\cdot g}& \C
}
\eeq
}\eit
\elm
\bpr
1. We should notice first that it is sufficient to prove this only for the case of $\mu=0$. Indeed, suppose we proved the following propositions:
\bit{
\item[$(i)^\circ$] the spaces $X_\lambda\cdot Y$ and $X\cdot Y$ become isomorphic after pseudosaturation:
\beq\label{(X_lambda-cdot-Y)^vartriangle=(X-cdot-Y)^vartriangle}
(X_\lambda\cdot Y)^\vartriangle=(X\cdot Y)^\vartriangle,
\eeq

\item[$(ii)^\circ$] for any functionals $f\in (X_\lambda)^\star$ and $g\in Y^\star$ the functional $f\cdot g:X_\lambda\cdot Y\to\C$ is continuous on the space $(X\cdot Y)^\vartriangle$, i.e. there exists a uniwue functional $h\in ((X\cdot Y)^\vartriangle)^\star$, such that the following diagram is commutative:
\beq\label{wedge_lambda^X-cdot-id_Y}
\xymatrix @R=2.pc @C=4.0pc 
{
(X_\lambda\cdot Y)^\vartriangle \ar[d]_{\vartriangle_{X_\lambda\cdot Y}} \ar[r]^{(\wedge_\lambda^X\cdot\id_Y)^\vartriangle} & (X\cdot Y)^\vartriangle \ar@{-->}[d]^{h} \\
X_\lambda\cdot Y \ar[r]_{f\cdot g}& \C
}
\eeq
}\eit
Then, first, \eqref{(X_lambda-cdot-Y_mu)^vartriangle=(X-cdot-Y)^vartriangle} holds, since
$$
(X_\lambda\cdot Y_\mu)^\vartriangle=\eqref{(X_lambda-cdot-Y)^vartriangle=(X-cdot-Y)^vartriangle}=(X\cdot Y_\mu)^\vartriangle=\eqref{X-cdot-Y-cong-Y-cdot-X}=
(Y_\mu\cdot X)^\vartriangle=\eqref{(X_lambda-cdot-Y)^vartriangle=(X-cdot-Y)^vartriangle}=
(Y\cdot X)^\vartriangle=\eqref{X-cdot-Y-cong-Y-cdot-X}=(X\cdot Y)^\vartriangle,
$$
and, second, for all functionals $f\in (X_\lambda)^\star$ and $g\in (Y_\mu)^\star$
the functional
$$
f\cdot g\ \circ\vartriangle_{X_\lambda\cdot Y_\mu}
$$
will be continuous on the space $(X_\lambda\cdot Y_\mu)^\vartriangle$, and therefore, on the (isomorphic to it by \eqref{(X_lambda-cdot-Y_mu)^vartriangle=(X-cdot-Y)^vartriangle}) space $(X\cdot Y)^\vartriangle$, and this functional on $(X\cdot Y)^\vartriangle$ is exactly the functional $h$, which is constructed in \eqref{wedge_lambda^X-cdot-wedge_mu^Y}. Its uniqueness follows also from \eqref{(X_lambda-cdot-Y_mu)^vartriangle=(X-cdot-Y)^vartriangle}.

2. Thus, we understood that it is sufficient to prove the weaker propositions $(i)^\circ$ and $(ii)^\circ$. They are proved by induction by ordinals  $\lambda\in\bf{Ord}$.

First, for $\lambda=0$ the equality \eqref{(X_lambda-cdot-Y_mu)^vartriangle=(X-cdot-Y)^vartriangle} becomes trivial
$$
(X\cdot Y)^\vartriangle=(X\cdot Y)^\vartriangle,
$$
and diagram \eqref{wedge_lambda^X-cdot-id_Y} turns into the diagram
$$
\xymatrix @R=2.pc @C=8.0pc 
{
(X\cdot Y)^\vartriangle \ar[d]_{\vartriangle_{X\cdot Y}} \ar[r]^{(\id_X\cdot\id_Y)^\vartriangle=\id_{X\cdot Y}^\vartriangle} & (X\cdot Y)^\vartriangle \ar@{-->}[d]^{h} \\
X\cdot Y \ar[r]_{f\cdot g}& \C
}
$$
where we just have to put $h=f\cdot g\circ\vartriangle_{X\cdot Y}$.

Further, suppose we proved $(i)^\circ$ and $(ii)^\circ$ for all ordinals $\iota<\lambda$, where $\lambda$ is a given ordinal, i.e.
\bit{
\item[$(i)^{\circ\circ}$] the spaces $X_\iota\cdot Y$ and $X\cdot Y$ become isomorphic after pseudosaturation:
\beq\label{(X_iota-cdot-Y)^vartriangle=(X-cdot-Y)^vartriangle}
(X_\iota\cdot Y)^\vartriangle=(X\cdot Y)^\vartriangle,
\eeq

\item[$(ii)^{\circ\circ}$] for any functionals $f\in (X_\iota)^\star$ and $g\in Y^\star$ the functional $f\cdot g:X_\iota\cdot Y\to\C$ is continuous on the space $(X\cdot Y)^\vartriangle$, i.e. there exists a unique functional $h\in ((X\cdot Y)^\vartriangle)^\star$, such that the following diagram is commutative:
\beq\label{wedge_iota^X-cdot-id_Y}
\xymatrix @R=2.pc @C=4.0pc 
{
(X_\iota\cdot Y)^\vartriangle \ar[d]_{\vartriangle_{X_\iota\cdot Y}} \ar[r]^{(\wedge_\iota^X\cdot\id_Y)^\vartriangle} & (X\cdot Y)^\vartriangle \ar@{-->}[d]^{h} \\
X_\iota\cdot Y\ar[r]_{f\cdot g}& \C
}
\eeq
}\eit
Let us show that $(i)^{\circ\circ}$ and $(ii)^{\circ\circ}$ hold for the ordinal $\lambda$ substituted insted of $\iota$ (in other words, we have to prove $(i)^{\circ}$ and $(ii)^{\circ}$). We have to consider two cases.

a) Suppose that $\lambda$ is an isolated ordinal, i.e.
$$
\lambda=\iota+1
$$
for some ordinal $\iota<\lambda$. Let us show first that $(ii)^\circ$ holds for $\iota=\lambda$. Take
$$
f\in (X_\lambda)^\star=(X_{\iota+1})^\star=((X_\iota)^\wedge)^\star=
\text{\cite[Theorem 3.10]{Ak03}}=
((X_\iota)^\star)^\vee
$$
and
$$
g\in Y^\star.
$$
Then $f$ is a limit of a totally bounded net $\{f_i;\ i\in I\}\subseteq (X_\iota)^\star$:
\beq\label{(X_lambda-cdot-Y_mu)^vartriangle=(X-cdot-Y)^vartriangle-1}
f_i\overset{((X_\iota)^\star)^\vee}{\underset{i\to\infty}{\longrightarrow}}f
\eeq
Since the set $\{f_i\}\subseteq (X_\iota)^\star$ is totally bouded, by Corollary \ref{LM:F-cdot-g-subseteq-(X-cdot-Y)^star}, the set
$$
\{f_i\cdot g;\ i\in I\}
$$
is totally bounded in $(X_\iota\cdot Y)^\star$. On the other hand, Theorem  \ref{TH:(f,g)|->f-cdot-g-razdelno-nepreryvno} implies that $\{f_i\cdot g;\ i\in I\}$ is a Cauchy net. Thus, we have that $\{f_i\cdot g;\ i\in I\}$ is a totally bounded Cauchy net in the space $(X_\iota\cdot Y)^\star$, hence, it converges in the enveloping space $(X_\iota\cdot Y)^{\star\triangledown}$:
$$
(X_\iota\cdot Y)^\star\subseteq (X_\iota\cdot Y)^{\star\triangledown}=
\text{\cite[Theorem 3.14]{Ak03}}=(X_\iota\cdot Y)^{\vartriangle\star}=\eqref{(X_iota-cdot-Y)^vartriangle=(X-cdot-Y)^vartriangle}=
(X\cdot Y)^{\vartriangle\star}
$$
In other words there is a functional $h\in (X\cdot Y)^{\vartriangle\star}$ such that
\beq\label{(X_lambda-cdot-Y_mu)^vartriangle=(X-cdot-Y)^vartriangle-2}
f_i\cdot g\overset{(X\cdot Y)^{\vartriangle\star}}{\underset{i\to\infty}{\longrightarrow}}h
\eeq
On the other hand, the chain of equalities
$$
((X_\iota)^\star)^\vee=\text{\cite[Theorem 3.10]{Ak03}}= ((X_\iota)^\wedge)^\star=(X_{\iota+1})^\star=
(X_\lambda)^\star
$$
implies that the relation \eqref{(X_lambda-cdot-Y_mu)^vartriangle=(X-cdot-Y)^vartriangle-1} is equivalent to the relation
\beq\label{(X_lambda-cdot-Y_mu)^vartriangle=(X-cdot-Y)^vartriangle-3}
f_i\overset{(X_\lambda)^\star}{\underset{i\to\infty}{\longrightarrow}}f
\eeq
which by Theorem \ref{TH:(f,g)|->f-cdot-g-razdelno-nepreryvno} implies
\beq\label{(X_lambda-cdot-Y_mu)^vartriangle=(X-cdot-Y)^vartriangle-4}
f_i\cdot g\overset{(X_\lambda\cdot Y)^{\star}}{\underset{i\to\infty}{\longrightarrow}} f\cdot g.
\eeq
Together \eqref{(X_lambda-cdot-Y_mu)^vartriangle=(X-cdot-Y)^vartriangle-2} and \eqref{(X_lambda-cdot-Y_mu)^vartriangle=(X-cdot-Y)^vartriangle-4} mean, in particular, that on each operator
$$
\ph\in X_\lambda\cdot Y=Y\underset{\xi}{:} (X_\iota)^\star=Y\underset{\xi}{:} X^\star=X_\lambda\cdot Y
$$
(by Lemma \ref{LM:o-drobyah-*}, the spaces $Y\underset{\xi}{:} (X_\iota)^\star$ and $Y\underset{\xi}{:} X^\star$ coincide as sets) the functionals $f\cdot g$ and $h$ coincide:
$$
(f\cdot g)(\ph)=h(\ph).
$$
This proves $(ii)^\circ$ (for the case of $\lambda=\iota+1$).

Let us prove now $(i)^\circ$. The mapping
$$
(X_\lambda\cdot Y)^\vartriangle=(Y\underset{\xi}{:}(X_\lambda)^\star)^\vartriangle\to (Y\underset{\xi}{:}X^\star)^\vartriangle=(X\cdot Y)^\vartriangle
$$
is always continuous. We have to prove that the inverse mapping
$$
(X\cdot Y)^\vartriangle=(Y\underset{\xi}{:}X^\star)^\vartriangle\to (Y\underset{\xi}{:}(X_\lambda)^\star)^\vartriangle=(X_\lambda\cdot Y)^\vartriangle
$$
is continuous. Since the space $(X\cdot Y)^\vartriangle$ is pseudosaturated, it is sufficient to prove the continuity of the mapping
$$
(X\cdot Y)^\vartriangle=(Y\underset{\xi}{:}X^\star)^\vartriangle\to Y\underset{\xi}{:}(X_\lambda)^\star=X_\lambda\cdot Y.
$$
Consider the basic neighbourhood of zero in $X_\lambda\cdot Y$, i.e. a set $U\cdot V$, where $U$ is a closed convex balanced neighbourhood of zero in $X_\lambda$, and $V$  is a closed convex balanced neighbourhood of zero in $Y$. By Lemma \ref{LM:o-drobyah}, the set $U\cdot V$ is capacious in $X\cdot Y$. Hence if we prove that it is closed in $(X\cdot Y)^\vartriangle$, this will mean that it is a neighbourhood of zero in $(X\cdot Y)^\vartriangle$.

This becomes obvious if we represent $U\cdot V=V:U^\circ$ as a polar of the system of functionals of the form $\{f\cdot g;\ f\in U^\circ,\ g\in V^\circ\}$:
\beq\label{(X_lambda-cdot-Y_mu)^vartriangle=(X-cdot-Y)^vartriangle-5}
U\cdot V=V:U^\circ=\{f\cdot g;\ f\in U^\circ,\ g\in V^\circ\}^\circ
\eeq
Indeed,
\begin{multline*}
\ph\in U\cdot V=V:U^\circ\quad\Leftrightarrow\quad \ph(U^\circ)\subseteq V
\quad\Leftrightarrow\quad \forall f\in U^\circ \ \ph(f)\in V
\quad\Leftrightarrow \\ \Leftrightarrow\quad
\forall g\in V^\circ\ \forall f\in U^\circ \ |\underbrace{g(\ph(f))}_{\scriptsize \begin{matrix}\|\\ (f\cdot g)(\ph)\end{matrix}}|\le 1
\quad\Leftrightarrow\quad \ph\in \{f\cdot g;\ f\in U^\circ,\ g\in V^\circ\}^\circ
\end{multline*}
Now we have to notice that in formula \eqref{(X_lambda-cdot-Y_mu)^vartriangle=(X-cdot-Y)^vartriangle-5}
$f\in U^\circ$, where $U$ is a neighbourhood of zero in $X_\lambda$, and $g\in V^\circ$, where $V$  is a neighbourhood of zero in $Y$. Hence, $f\in (X_\lambda)^\star$, and $g\in Y^\star$. Thus, by $(ii)^\circ$ (which is already proved for the case of $\lambda=\iota+1$), $f\cdot g$ are continuous functionals on the space $(X\cdot Y)^\vartriangle$. We see that $U\cdot V$ is a polar of some system of continuous functionals on $(X\cdot Y)^\vartriangle$, hence, $U\cdot V$ is a closed set in $(X\cdot Y)^\vartriangle$. In addition it is capacious, hence it is a neighbourhood of zero in the (pseudosaturated) space $(X\cdot Y)^\vartriangle$. We proved $(i)^\circ$ (for the case of $\lambda=\iota+1$).

b) Now we consider the case when $\lambda$ is a limit ordinal, i.e.
$$
\lambda\ne \iota+1
$$
for all $\iota<\lambda$. Again we have to prove first $(ii)^\circ$. Take $f\in (X_\lambda)^\star$ and $g\in Y^\star$. Then
$$
f\in (X_\lambda)^\star=\text{\cite[Theorem 3.12]{Ak03}}
=(X^\star)^\lambda= \text{\cite[(1.32)]{Ak03}}=\bigcup_{\iota<\lambda}(X^\star)^\iota,
$$
hence there is an ordinal $\iota<\lambda$ such that $f\in (X^\star)^\iota=\text{\cite[Theorem 3.12]{Ak03}}=(X_\iota)^\star$. We see that $f\in (X_\iota)^\star$ and $g\in Y^\star$, hence by the assumpotion of the induction in $(ii)^{\circ\circ}$, the functional $f\cdot g$ must be continuous on the space $(X\cdot Y)^\vartriangle$. This proves $(ii)^\circ$ (for the case of the limit $\lambda$).

Now we can prove $(i)^\circ$. We repeat here the reasoning of a): we have to prove the continuity of the mapping
$$
(X\cdot Y)^\vartriangle=(Y\underset{\xi}{:}X^\star)^\vartriangle\to Y\underset{\xi}{:}(X_\lambda)^\star=X_\lambda\cdot Y
$$
We take a basic neighbourhood of zero in $X_\lambda\cdot Y$, i.e. a set $U\cdot V$, where $U$ is a neighbourhood of zero in $X_\lambda$, and $V$ a neighbourhood of zero in $Y$, then we prove equality
\eqref{(X_lambda-cdot-Y_mu)^vartriangle=(X-cdot-Y)^vartriangle-5}, then we notice that in it $f\in (X_\lambda)^\star$, and $g\in Y^\star$. After that by already proven $(ii)^\circ$ (for the case of the limit $\lambda$), we have that $f\cdot g$ are continuous functionals on the space $(X\cdot Y)^\vartriangle$. Hence, $U\cdot V$ is a closed set in $(X\cdot Y)^\vartriangle$. In addition, it is capacious, by Lemma \ref{LM:o-drobyah}, hence it is a neighbourhood of zero in the (pseudosaturated) space $(X\cdot Y)^\vartriangle$. We have proved $(i)^\circ$ (for the case of the limit $\lambda$).
\epr

For arbitrary locally convex spaces $X$ and $Y$ the mappings
$$
\vartriangle_X: X^\vartriangle\to X,\qquad \vartriangle_Y: Y^\vartriangle\to Y
$$
define the mapping
\beq\label{vartriangle_X-cdot-vartriangle_Y}
(\vartriangle_X\cdot \vartriangle_Y): (X^\vartriangle\cdot Y^\vartriangle) \to (X\cdot Y)
\eeq
(which coincide with the mapping \eqref{wedge_lambda^X-cdot-wedge_mu^Y-1} for some $\lambda$ and $\mu$).

\btm\label{TH:(X^vartriangle-cdot-Y^vartriangle)^vartriangle=(X-cdot-Y)^-vartriangle}
For each pseudocomplete locally convex spaces $X$ and $Y$ the pseudosaturation of the mapping \eqref{vartriangle_X-cdot-vartriangle_Y} is an isomorphism of locally convex (and stereotype) spaces:
\beq\label{(X^vartriangle-cdot-Y^vartriangle)^vartriangle=(X-cdot-Y)^-vartriangle}
X^\vartriangle\odot Y^\vartriangle\cong (X^\vartriangle\cdot Y^\vartriangle)^\vartriangle\cong (X\cdot Y)^\vartriangle.
\eeq
\etm
\bpr
The second isomorphism follows from \eqref{(X_lambda-cdot-Y_mu)^vartriangle=(X-cdot-Y)^vartriangle}, if we choose ordinals $\lambda$ and $\mu$ such that $X_\lambda=X^\vartriangle$ and $Y_\mu=Y^\vartriangle$. And the first isomorphism uses the fact that the pseudosaturation $X^\vartriangle$ of a pseudocompolete space $X$ is always pseudocomplete (and thus, stereotype) \cite[Proposition 3.16]{Ak03}. From this we have that in the space $(X^\vartriangle)^\star$ each totally bounded set $S$ is contained in the polar $U^\circ$ of some neighbourhood of zero $U\subseteq X^\vartriangle$. As a corollary, in the space of operators $Y^\vartriangle:(X^\vartriangle)^\star$ the topology of uniform convergence on totally bounded sets coincides with the topology of uniform convergence on polars $U^\circ$ of neighbourhoods of zero $U\subseteq X^\vartriangle$:
$$
Y^\vartriangle:(X^\vartriangle)^\star=Y^\vartriangle\underset{\xi}{:}(X^\vartriangle)^\star
$$
Thin implies the following chain that proves the first isomorphism in \eqref{(X^vartriangle-cdot-Y^vartriangle)^vartriangle=(X-cdot-Y)^-vartriangle}:
$$
X^\vartriangle\odot Y^\vartriangle=Y^\vartriangle\oslash (X^\vartriangle)^\star=
(Y^\vartriangle:(X^\vartriangle)^\star)^\vartriangle=
(Y^\vartriangle\underset{\xi}{:}(X^\vartriangle)^\star)^\vartriangle= \eqref{DEF:X-cdot-Y}=(X^\vartriangle\cdot Y^\vartriangle)^\vartriangle.
$$
\epr

\tableofcontents

\end{document}


\title{О тензорных дробях и тензорных произведениях в категории стереотипных пространств}

\author{С.С.Акбаров}

\address{Департамент прикладной математики, Национальный исследовательский институт Высшая школа экономики, Таллинская 34, Москва, 123458 Россия}
\email{sergei.akbarov@gmail.com}
\thanks{Поддержано грантом РФФИ № 18-01-00398.}

\maketitle

\begin{abstract}
В статье доказываются два тождества, связывающие некоторые естественные тензорные произведения в категории $\sf{LCS}$ локально выпуклых пространств с тензорными произведениями в категории $\sf{Ste}$ стереотипных пространств. В частности, приводятся условия, при которых выполняется тождество
$$
X^\vartriangle\odot Y^\vartriangle\cong (X^\vartriangle\cdot Y^\vartriangle)^\vartriangle\cong (X\cdot Y)^\vartriangle,
$$
в котором $\odot$ --- инъективное тензорное произведение в категории $\sf{Ste}$, $\cdot$ --- первичное тензорное произведение в категории $\sf{LCS}$, а $\vartriangle$ --- операция псевдонасыщения в категории $\sf{LCS}$. Изучение соотношений этого типа оправдывается тем, что они оказываются важными инструментами при построении теорий двойственности, основанных на понятии оболочки. В частности, с их помощью строится обобщение голоморфной теории двойственности на класс (необязательно абелевых) счетных дискретных групп.
\end{abstract}

\section{Введение}

В работах автора \cite{Ak95-1,Ak95-2,Ak03,Ak08,Ak16} были описаны свойства класса $\Ste$ локально выпуклых пространств, названных {\it стереотипными}  и определяемых условием\footnote{В книге Г.~Кете \cite{Kothe} пространства такого типа называются {\it полярно рефлексивными}.}
$$
X\cong (X^\star)^\star,
$$
в котором каждая звездочка $\star$ означает сопряженное пространство функционалов, наделенное топологией равномерной сходимости на вполне ограниченных множествах. Как отмечалось, этот класс обладает серией замечательных свойств, в частности,
\bit{

\item[---] $\Ste$ весьма широк, потому что содержит все квазиполные бочечные  пространства (в частности, все пространства Фреше и поэтому все банаховы пространства),

\item[---] $\Ste$ образует предабелеву, полную и кополную категорию с линейными непрерывными отображениями в качестве морфизмов,

\item[---] категория $\Ste$ автодуальна относительно функтора $X\mapsto X^\star$,

\item[---] $\Ste$ обладает тремя естественными бифункторами $\circledast$, $\odot$, $\oslash$ со следующим свойствами:
\begin{eqnarray}
 \C\circledast X \cong X  \cong  X\circledast \C, & \C\odot X \cong  X  \cong  X\odot \C, \label{C-circledast-X-cong-X}
\\
 X\circledast Y \cong  Y\circledast X, & X\odot Y \cong Y\odot
X, \label{X-circledast-Y-cong-Y-circledast-X}
\\
 (X\circledast Y)\circledast Z  \cong  X\circledast
(Y\circledast Z), &  (X\odot Y)\odot Z  \cong  X\odot (Y\odot Z),
\label{X-circledast-Y-circledast-Z}
 \\
    Z\oslash (Y\circledast X)\cong (Z\oslash Y)\oslash X,
 &  (X\odot Y)\oslash Z \cong X\odot (Y\oslash Z), \label{Z-oslash-(Y-circledast-X)} \\
(X\circledast Y)^\star  \cong   Y^\star \odot X^\star, &
(X\odot Y)^\star \cong   Y^\star \circledast X^\star, \\
 X\circledast Y\cong (X^\star\oslash Y)^\star, & X\odot Y\cong Y\oslash (X^\star).
\end{eqnarray}
}\eit
Тождества \eqref{C-circledast-X-cong-X}---\eqref{X-circledast-Y-circledast-Z} в этом списке означают, что $\Ste$ является симметрической моноидальной категорией относительно бифункторов $\circledast$ и $\odot$, и это оправдывает для них название {\it ``стереотипные тензорные произведения''} ($\circledast$ --- {\it проективное}, а $\odot$ --- {\it инъективное}) . С другой стороны, тождества \eqref{Z-oslash-(Y-circledast-X)} оправдывают для бифунктора $\oslash$ название {\it ``стереотипная тензорная дробь''} (при этом левое из них означает, что относительно $\circledast$ моноидальная категория $\Ste$ замкнута с внутренним hom-функтором $\oslash$).

В дальнейших исследованиях автором (а позднее Ю.Н.Кузнецовой и О.Ю.Аристовым) отмечались многочисленные приложения категории стереотипных пространств в функциональном анализе и геометрии, и, в частности, в построении обобщений понтрягинской двойственности на некоммутативные группы. Как один из примеров, из результатов работ \cite{Ak08,ArHRC} следует, что понтрягинская двойственность продолжается с категории конечных абелевых групп на категорию аффинных алгебраических групп функтором $G\mapsto {\mathcal O}^\star(G)$ перехода к групповой алгебре аналитических функционалов. Это обобщение выражается функториальной диаграммой:
{\sf
 \beq\label{diagramma-kategorij-dlya-konechnyh-affinnyh-grupp}
 \xymatrix @R=3.pc @C=2.pc
 {
 \boxed{\begin{matrix}
  \text{алгебры Хопфа,}\\
  \text{рефлексивные относительно}\\
  \text{оболочки Аренса---Майкла}
 \end{matrix}}
 \ar[rr]^{H\mapsto H^\dagger} & &
 \boxed{\begin{matrix}
  \text{алгебры Хопфа,}\\
  \text{рефлексивные относительно}\\
  \text{оболочки Аренса---Майкла}
 \end{matrix}}
 \\
 \boxed{\begin{matrix}
  \text{аффинные алгебраические группы}
 \end{matrix}} \ar[u]^(0.4){\scriptsize\begin{matrix} {\mathcal O}^\star(G)\\
 \text{\rotatebox{90}{$\mapsto$}} \\ G\end{matrix}} & &
 \boxed{\begin{matrix}
  \text{аффинные алгебраические группы}
 \end{matrix}} \ar[u]_(0.4){\scriptsize\begin{matrix} {\mathcal O}^\star(G) \\
 \text{\rotatebox{90}{$\mapsto$}} \\ G\end{matrix}} \\
 \boxed{\begin{matrix}
 \text{абелевы конечные группы}
 \end{matrix}} \ar[u]^{\mathfrak{e}} \ar[rr]^{G\mapsto G^\bullet} & &
  \boxed{\begin{matrix}
 \text{абелевы конечные группы}
  \end{matrix}}\ar[u]_{\mathfrak{e}}
 }
 \eeq }\noindent
в которой $\mathfrak{e}$ обозначает естественное вложение категорий, $G^\bullet$ --- двойственную по Понтрягину группу, $\dagger$ --- композицию функтора $H\mapsto \widehat{H}$ оболочки Аренса---Майкла (см. определение в \cite{Ak16} или в \cite{AkCountGroups}) и функтора $X\mapsto X^\star$ перехода к сопряженному стереотипному пространству,
$$
H^\dagger:=( \widehat{H} )^\star,
$$
а под алгебрами Хопфа, рефлексивными относительно оболочки Аренса---Майкла, понимаются алгебры Хопфа в моноидальной категории $(\Ste,\circledast)$, удовлетворяющие тождеству
$$
H\cong \big(H^\dagger\big)^\dagger
$$
(см. подробности в \cite{AkCountGroups}). Коммутативность же диаграммы \eqref{diagramma-kategorij-dlya-konechnyh-affinnyh-grupp} понимается как изоморфизм двух функторов, ведущих из левого нижнего угла в правый верхний:
\beq\label{O^star(G)^dagger-cong-O^star(G^bullet)}
{\mathcal O}^\star(G)^\dagger\cong {\mathcal O}^\star(G^\bullet)
\eeq

После работы автора \cite{Ak08}, в которой двойственность относительно  оболочки была впервые описана и применена к комплексным группам Ли, эта же идея (с заменой оболочки Аренса---Майкла на другие оболочки) была использована для обобщений понтрягинской двойственности на разные другие классы групп. В частности, Ю.Н.Кузнецова в работе \cite{Kuz} предложила теорию двойственности для (необязательно коммутативных) групп Мура (впоследствии ее результаты уточнялись в \cite{Ak17-1,Ak17-2,Akbarov-Moore}), и, кроме того, сам автор в \cite{Ak17-1,Ak17-2} описал теорию двойственности для дифференциальной геометрии, работающую по крайней мере для вещественных групп Ли вида $\R^n\times K\times D$, где $K$ --- компактная группа Ли, а $D$ --- дискретная группа Мура.

Эти наблюдения не получили пока оформления в виде теорий с естественными границами, потому что очевидно, что все эти результаты можно обобщать на более широкие классы групп, и неясно, где в этих исследованиях логично будет остановиться. В рамках этой программы в работе \cite{AkCountGroups} автором был рассмотрен вопрос, можно ли по аналогии с диаграммой \eqref{diagramma-kategorij-dlya-konechnyh-affinnyh-grupp} обобщить понтрягинскую двойственность с категории конечных абелевых групп на категорию счетных дискретных групп. В \cite{AkCountGroups} на него был дан положительный ответ, описывающийся диаграммой
{\sf
 \beq\label{diagramma-kategorij-dlya-diskretnyh-grupp}
 \xymatrix @R=3.pc @C=2.pc
 {
 \boxed{ \begin{matrix}
  \text{голоморфно рефлексивные}\\
  \text{алгебры Хопфа}\\
 \end{matrix}}
 \ar[rr]^{H\mapsto H^\dagger} & &
 \boxed{ \begin{matrix}
  \text{голоморфно рефлексивные}\\
  \text{алгебры Хопфа}\\
 \end{matrix}}
 \\
 \boxed{\begin{matrix}
 \text{счетные дискретные группы}
 \end{matrix}} \ar[u]^(0.4){\scriptsize\begin{matrix} {\mathcal O}^\star(G)\\
 \text{\rotatebox{90}{$\mapsto$}} \\ G\end{matrix}} & &
 \boxed{\begin{matrix}
 \text{счетные дискретные группы}
 \end{matrix}} \ar[u]_(0.4){\scriptsize\begin{matrix} {\mathcal O}^\star(G) \\
 \text{\rotatebox{90}{$\mapsto$}} \\ G\end{matrix}} \\
 \boxed{\begin{matrix}
 \text{конечные абелевы группы}
 \end{matrix}} \ar[u]^{\mathfrak{e}} \ar[rr]^{G\mapsto G^\bullet} & &
  \boxed{\begin{matrix}
 \text{конечные абелевы группы}
  \end{matrix}}\ar[u]_{\mathfrak{e}}
 }
 \eeq }\noindent
в которой под голоморфной рефлексивностью алгебры Хопфа понимается ее рефлексивность  относительно одного из стереотипных вариантов оболочки Аренса--Майкла (см. детали в \cite{AkCountGroups}).

Доказательство этого факта, однако, как оказалось, требует глубокого погружения в основания теории стереотипных пространств с выводом некоей серии технических утверждений, которые из-за их объема и идеологической обособленности, желательно выделить в отдельный текст. Настоящая статья как раз является таким текстом. Здесь мы выводим несколько формул, выражающих инъективное стереотипное тензорное произведение $\odot$ и тензорную дробь $\oslash$, упоминавшиеся выше, через некоторые более простые конструкции в теории топологических векторных пространств.

Первая формула выражает $\odot$ через один из вариантов тензорного произведения в категории локально выпуклых пространств, который мы обозначаем точкой $\cdot$ и называем ``первичным тензорным произведением''. Формально $X\cdot Y$ определяется как пространство линейных непрерывных отображений $\ph:X^\star\to Y$ с топологией равномерной сходимости на полярах $U^\circ$ окрестностей нуля в $X$, и для него мы доказываем формулу
\beq\label{(X^vartriangle-cdot-Y^vartriangle)^vartriangle=(X-cdot-Y)^-vartriangle-0}
X^\vartriangle\odot Y^\vartriangle\cong (X^\vartriangle\cdot Y^\vartriangle)^\vartriangle\cong (X\cdot Y)^\vartriangle.
\eeq
справедливую для произвольных псевдополных локально выпуклых пространств $X$ и $Y$ (без условия, что какое-то из них должно обладать классической аппроксимацией).

Вторая же формула дополняет формальное определение тензорной дроби $\oslash$ через пространство операторов. Напомним, что в \cite{Ak03} пространство $\oslash$ определялось в два этапа. Сначала вводилось пространство $Y:X$ (``первичная тензорная дробь''), под которым понималось пространство линейных непрерывных отображений $\ph:X\to Y$ с топологией равномерной сходимости на вполне ограниченных множествах. После этого пространство $Y\oslash X$ определяется как псевдонасыщение пространства $Y:X$
\beq\label{DEF:oslash}
    Y\oslash X=(Y:X)^\vartriangle.
\eeq
Так вот формула, связывающая $\oslash$ с $:$, и доказываемая в этой статье, выглядит так:
\beq\label{(vartriangle_Y:triangledown_X)^vartriangle-isomorphism-0}
Y^\vartriangle\oslash X^\triangledown\cong(Y^\vartriangle:X^\triangledown)^\vartriangle\cong  (Y:X)^\vartriangle
\eeq
где $\triangledown$ --- операция псевдопополнения в категории локально выпуклых пространств \cite[1.3]{Ak03} (а $\vartriangle$ --- по-прежнему, операция псевдонасыщения).

Формула
\eqref{(X^vartriangle-cdot-Y^vartriangle)^vartriangle=(X-cdot-Y)^-vartriangle-0} ценна тем, что позволяют судить о свойствах пространства $X^\vartriangle\odot Y^\vartriangle$, без детального исследования свойств топологий псевдонасыщений $X^\vartriangle$ и $Y^\vartriangle$. Например, в \cite{AkCountGroups} эти формулы применяются для доказательства тождеств вида
\beq\label{F(M)-odot-F(N)}
{\mathcal F}(M)^\vartriangle\odot {\mathcal F}(N)^\vartriangle\cong {\mathcal F}(M\times N)^\vartriangle
\eeq
в которых $M$ и $N$ --- многообразия определенного типа, а ${\mathcal F}(...)$ --- определенные типы (топологических векторных) функциональных пространств на таких многообразиях. Типичная ситуация в теории стереотипных пространств --- когда мы легко можем описать топологию какого-то пространства, в данном случае, ${\mathcal F}(...)$, но топология его псевдонасыщения ${\mathcal F}(...)^\vartriangle$ выглядит необычайно сложной, настолько, что ничего проще самого определения пседонасыщения для ее описания найти обычно невозможно. Так вот, формула \eqref{F(M)-odot-F(N)} дает описание стереотипного тензорного произведения без необходимости анализировать, как устроены топологии пространств ${\mathcal F}(...)^\vartriangle$.

Точно так же формула \eqref{(vartriangle_Y:triangledown_X)^vartriangle-isomorphism-0} позволяет свести описание пространства $Y^\vartriangle\oslash X^\triangledown$ к описанию свойств пространства $Y:X$ без необходимости изучать свойства пространств $Y^\vartriangle$ и $X^\triangledown$ (формулу \eqref{(vartriangle_Y:triangledown_X)^vartriangle-isomorphism-0}, мы, впрочем, приводим как попутный результат, потому что напрямую в \cite{AkCountGroups} она не используется).

Этим вопросам посвящена настоящая работа.

\subsection{Благодарности}

Автор благодарит О.Ю.Аристова за полезные обсуждения.

\section{Псевдонасыщение первичной дроби $Y:X$}

Всюду в тексте мы используем терминологию и обозначения работы \cite{Ak03}. Мы рассматриваем локально выпуклые пространства над полем $\C$. Для всякого локально выпуклого пространства $X$ символ ${\mathcal U}(X)$ обозначает множество всех окрестностей нуля в $X$, ${\mathcal S}(X)$ --- множество всех вполне ограниченных подмножеств в $X$, ${\mathcal D}(X)$ --- множество всех емких\footnote{Множество $D\subseteq X$ называется {\it емким}, если для любого вполне ограниченного множества $S\subseteq X$ найдется конечное множество $F\subseteq X$, такое что $S\subseteq D+F$.} подмножеств в $X$. Термин {\it оператор} будет использоваться для линейных непрерывных отображений $\ph:X\to Y$ локально выпуклых пространств.

Мы начнем с доказательства формулы \eqref{(vartriangle_Y:triangledown_X)^vartriangle-isomorphism-0}, последней из анонсированных во Введении.

\subsection{Билинейная форма $(g,x)\mapsto g:x$.}

Напомним, что в \cite[5.6]{Ak03} билинейное отображение $\beta:X\times Y\to Z$ называлось {\it непрерывным}, если
\bit{
\item[1)] для всякой окрестности нуля $W\subseteq Z$ и любого вполне ограниченного множества $S\subseteq X$ найдется окрестность нуля $V\subseteq Y$ такая, что
    $$
    \beta(S,V)\subseteq W,
    $$
    и
\item[2)] для всякой окрестности нуля $W\subseteq Z$ и любого вполне ограниченного множества $T\subseteq Y$ найдется окрестность нуля $U\subseteq X$ такая, что
    $$
    \beta(U,T)\subseteq W.
    $$

}\eit

Еще раз напомним, что для любых локально выпуклых пространств $X$ и $Y$ символ $Y:X$ ({\it ``первичная тензорная дробь''}) обозначает пространство линейных непрерывных отображений $\ph:X\to Y$ с топологией {\it равномерной сходимости на вполне ограниченных множествах} \cite[5.1]{Ak03}. Если $A\subseteq X$ и $B\subseteq Y$, то символ $B:A$ обозначает позмножество в $Y:X$, состоящее из отображений, переводящих $A$ в $B$ \cite[5.4]{Ak03}:
\beq\label{DEF:B:A}
\ph\in B:A\quad\Leftrightarrow\quad \ph\in Y:X\ \& \ \ph(A)\subseteq B.
\eeq
Для любых $x\in X$ и $g\in Y^\star$ определим функционал
\beq\label{DEF:(g:x)}
    (g:x):(Y:X)\to\C\quad\Big|\quad (g:x)(\ph)=g(\ph(x)),\qquad \ph\in Y:X
\eeq

\btm\label{TH:(g,x)|->g:x-razdelno-nepreryvno}
Пусть $X$ и $Y$ --- локально выпуклые пространства, причем $X$ псевдонасыщено. Тогда билинейное отображение
$$
(g,x)\in Y^\star\times X\mapsto g:x\in (Y:X)^\star
$$
непрерывно (в смысле определения выше).
\etm
\bpr
1. Пусть $g_i\to 0$ и $S\in{\mathcal S}(X)$. Тогда для произвольного вполне ограниченного множества $\varPhi\subseteq Y:X$ множество $\varPhi(S)$ будет в силу \cite[Theorem 5.1]{Ak03}  вполне ограничено в $Y$. Поэтому направленность $g_i$ будет стремиться к нулю равномерно на множестве $\varPhi(S)$ (по определению топологии в $Y^\star$):
$$
(g_i:x)(\ph)=g_i(\ph(x)) \overset{\C}{\underset{i\to\infty}{\underset{\ph\in\varPhi,\ x\in S}{\rightrightarrows}}}0.
$$
Это верно для любого вполне ограниченного множества $\varPhi\subseteq Y:X$, поэтому
$$
g_i:x \overset{Y:X}{\underset{i\to\infty}{\underset{x\in S}{\rightrightarrows}}}0.
$$

2. Пусть $G\in {\mathcal S}(Y^\star)$ и $x_i\to 0$. Тогда для произвольного вполне ограниченного множества $\varPhi\subseteq Y:X$ множество $G\circ\varPhi$ будет в силу \cite[Theorem 5.1]{Ak03}  вполне ограничено в $X^\star$. Поскольку по условию, $X$ псевдонасыщено, множество $G\circ\varPhi$  равностепенно непрерывно на $X$. Как следствие,
$$
(g:x_i)(\ph)=g(\ph(x_i)) \overset{\C}{\underset{i\to\infty}{\underset{\ph\in\varPhi,\ g\in G}{\rightrightarrows}}}0.
$$
Это верно для любого вполне ограниченного множества $\varPhi\subseteq Y:X$, поэтому
$$
g:x_i \overset{Y:X}{\underset{i\to\infty}{\underset{g\in G}{\rightrightarrows}}}0.
$$
\epr

\bcor\label{COR:(g,x)|->g:x-razdelno-nepreryvno}
Пусть $X$ и $Y$ --- локально выпуклые пространства, причем $X$ псевдонасыщено. Тогда если $S\subseteq X$ --- вполне ограниченное множество и $g\in Y^\star$, то  $g:S\subseteq (Y:X)^\star$ --- вполне ограниченное множество. Двойственным образом, если $x\in X$ и $G\subseteq Y^\star$ --- вполне ограниченное множество, то $G:x\subseteq (Y:X)^\star$ --- вполне ограниченное множество.
\ecor

\subsection{Псевдонасыщение $Y:X$.}

Пусть $\{X^\lambda;\ \lambda\in\bf{Ord}\}$ -- инъективный ряд пространства $X$, а $\{Y_\mu;\ \mu\in\bf{Ord}\}$ --- проективный ряд пространства $Y$, и
$$
\vee_\lambda^X: X\to X^\lambda,\qquad \wedge_\mu^Y: Y_\mu\to Y
$$
--- соответствующие естественные отображения. Рассмотрев аналог диаграммы \cite[(5.1)]{Ak03}
\beq\label{wedge_mu^Y:vee_lambda^X}
\xymatrix @R=2.pc @C=2.0pc 
{
X \ar[d]_{(\wedge_\mu^Y:\vee_\lambda^X)(\psi)= \wedge_\mu^Y\circ\psi\circ \vee_\lambda^X} \ar[r]^{\vee_\lambda^X} & X^\lambda \ar[d]^{\psi} \\
Y & Y_\mu \ar[l]_{\wedge_\mu^Y}
}
\eeq
мы получим непрерывное отображение
\beq\label{wedge_mu^Y:vee_lambda^X-1}
(\wedge_\mu^Y:\vee_\lambda^X): (Y_\mu:X^\lambda) \to (Y:X).
\eeq

\blm\label{LM:Y_mu:X^lambda=Y:X}
Если $X$ псевдонасыщено, а $Y$ псевдополно, то отображение \eqref{wedge_mu^Y:vee_lambda^X-1} является биекцией, поэтому пространства $Y_\mu:X^\lambda$ и $Y:X$ можно отождествить как множества (а как локально выпуклые пространства они будут отличаться только топологией):
\beq\label{Y_mu:X^lambda=Y:X}
Y_\mu:X^\lambda=Y:X
\eeq
\elm
\bpr
Отображение $\vee_\lambda^X: X\to X^\lambda$ является эпиморфизмом. С другой стороны, отображение $\wedge_\mu^Y: Y_\mu\to Y$ является инъекцией. Отсюда следует, что отображение
$$
\psi\mapsto \wedge_\mu^Y\circ\psi\circ (\wedge_\lambda^X)^\star
$$
является инъекцией. Нам нужно проверить, что оно является сюръекцией. Действительно, пусть $\ph:X\to Y$ --- произвольный оператор. Мы последовательно построим два оператора $\ph_\mu$ и $\psi$, замыкающие диаграмму:
$$
\xymatrix @R=3.pc @C=4.0pc 
{
X \ar[d]_{\ph}\ar@{-->}[dr]_{\ph_\mu} \ar[r]^{\vee_\lambda^X} & X^\lambda \ar@{-->}[d]^{\psi} \\
Y & Y_\mu \ar[l]^{\wedge_\mu^Y}
}
$$
Поскольку $X$ псевдонасыщено, оператор $\ph$ однозначно продолжается до оператора $\ph^\vartriangle:X\to Y^\vartriangle$. Если теперь рассмотреть его композицию с естественным отображением $Y^\vartriangle\to Y_\mu$, то мы получим оператор, продолжающий $\ph$ до оператора $\ph_\mu:X\to Y_\mu$.

Далее, поскольку $Y$ псевдополно, $Y_\mu$ тоже псевдополно (в силу \cite[Proposition 3.16]{Ak03}). Поэтому оператор  $\ph_\mu:X\to Y_\mu$ продолжается до некоторого оператора $\chi:X^\triangledown\to Y_\mu$. Рассмотрев его композицию с естественным вложением $X^\lambda\to X^\triangledown$, мы получим оператор $\psi:X^\lambda\to Y_\mu$.
\epr

\blm\label{LM:o-sohranenii-polnoi-ogranichennosti-v-Y_mu:X_lambda}
Пусть $X$ псевдонасыщено, а $Y$ псевдополно. Если $\varPhi\in Y:X$ --- вполне ограниченное множество, то для любых $\lambda,\mu\in\bf{Ord}$ его представление в пространстве $Y_\mu:X^\lambda$ (биекцией \eqref{Y_mu:X^lambda=Y:X}) --- также вполне ограниченное множество.
\elm
\bpr
1. Рассмотрим сначала случай, когда $\lambda=0$. Если $\varPhi\subseteq Y:X$ --- вполне ограниченное множество, то в силу \cite[Theorem 5.1]{Ak03} это значит, что  $\varPhi$ равностепенно непрерывно и равномерно вполне ограничено на каждом вполне ограниченном множестве $S\subseteq X$. Как следствие, образ $\varPhi(S)$ будет вполне ограниченным множеством в $Y$. Пространство $Y_\mu$ представляет собой усиление топологии на $Y$, при котором класс вполне ограниченных множеств и топология на вполне ограниченных множествах не меняются. Поэтому $\varPhi(S)$ будет вполне ограниченным множеством также и в пространстве $Y_\mu$, причем с той же топологией, что и индуцированная из $Y$. Отсюда можно заключить, что множество $\varPhi$, рассматриваемое как множество операторов из $X$ в $Y_\mu$, равностепенно непрерывно и равномерно вполне ограничено на $S$. И это верно для всякого вполне ограниченного множества $S\subseteq X$. Значит (опять, в силу \cite[Theorem 5.1]{Ak03}), $\varPhi$ вполне ограничено в $Y_\mu:X$.

2. Итак, мы поняли, что если $\varPhi\subseteq Y:X$ вполне ограничено в $Y:X$, то оно вполне ограничено и в $Y_\mu:X$ при любом $\mu\in\bf{Ord}$. Если взять достаточной большой ординал $\mu\in\bf{Ord}$, то мы получим, что $\varPhi\subseteq Y:X$ вполне ограничено в пространстве $Y^\vartriangle:X$. Отсюда, в силу \cite[Lemma 5.10]{Ak03}, мы получаем, что $\varPhi\subseteq Y:X$ вполне ограничено в пространстве $Y^{\vartriangle\triangledown}:X^\triangledown$.

Поскольку, по условию, $Y$ псевдополно, в силу \cite[Proposition 3.16]{Ak03} пространство $Y^{\vartriangle}$ тоже псевдополно, и, как следствие,
$$
Y^{\vartriangle\triangledown}=Y^{\vartriangle}.
$$
Мы получаем, что множество $\varPhi\subseteq Y:X$ вполне ограничено в пространстве $Y^\vartriangle:X^\triangledown$.

Теперь, зафиксировав произвольные $\lambda,\mu\in{\bf Ord}$ и обозначив через
$$
\sigma:X^\lambda\to X^\triangledown,\qquad
\pi:Y_\mu\gets Y^\vartriangle,
$$
естественные отображения, мы получим линейное непрерывное отображение
$$
(\pi:\sigma): (Y^\vartriangle:X^\triangledown) \to (Y_\mu:X^\lambda).
$$
Оно переводит вполне ограниченное множество $\varPhi\subseteq Y^\vartriangle:X^\triangledown$ во вполне ограниченное множество
$\varPhi\subseteq Y_\mu:X^\lambda$.
\epr

\blm\label{LM:o-drobyah-v-Y_mu:X_lambda}
Пусть $X$ псевдонасыщено, а $Y$ псевдополно. Тогда если $S\in {\mathcal S}(X^\lambda)$ и $V\in {\mathcal U}(Y_\mu)$, то $V:S\in {\mathcal D}(Y:X)$.\footnote{Обозначение $V:S$ введено в \eqref{DEF:B:A}.}
\elm
\bpr
Пусть $\varPhi\subseteq Y:X$ --- вполне ограниченное множество. По лемме \ref{LM:o-sohranenii-polnoi-ogranichennosti-v-Y_mu:X_lambda}, $\varPhi$ будет вполне ограничено и в пространстве $Y_\mu:X^\lambda$. Поскольку $V:S$ --- окрестность нуля в пространстве $Y_\mu:X^\lambda$, найдется конечное множество $A\subseteq Y_\mu:X^\lambda$ такое что
\beq\label{o-drobyah-v-Y_mu:X_lambda}
\varPhi\subseteq V:S+A.
\eeq
Но по лемме \ref{LM:Y_mu:X^lambda=Y:X}, все эти множества можно считать содержащимися в пространстве $Y:X$. Мы получем такую вещь: для любого вполне ограниченного множества $\varPhi$ в $Y:X$ существует конечное множество $A$ в $Y:X$ такое, что выполняется \eqref{o-drobyah-v-Y_mu:X_lambda}. Это значит, что множество $V:S$ является емким в пространстве $Y:X$.
\epr

\blm\label{LM:(Y_mu:X_lambda)^vartriangle=(Y:X)^vartriangle}
Если $X$ псевдонасыщено, а $Y$ псевдополно, то для любых ординалов $\lambda,\mu\in\bf{Ord}$
\bit{
\item[(i)] при псевдонасыщении пространства $Y_\mu:X^\lambda$ и $Y:X$ становятся изоморфными:
\beq\label{(Y_mu:X_lambda)^vartriangle=(Y:X)^vartriangle}
(Y_\mu:X^\lambda)^\vartriangle=(Y:X)^\vartriangle,
\eeq

\item[(ii)] для любых $x\in X^\lambda$ и $g\in (Y_\mu)^\star$ функционал
$$
    (g:x):(Y_\mu:X^\lambda)\to\C\quad\Big|\quad (g:x)(\ph)=g(\ph(x)),\qquad \ph\in Y_\mu:X^\lambda
$$
    непрерывен на пространстве $(Y:X)^\vartriangle$, то есть  существует единственный функционал $h\in ((Y:X)^\vartriangle)^\star$, замыкающий диаграмму
\beq\label{(Y_mu:X_lambda)^vartriangle=(Y:X)^vartriangle-func}
\xymatrix @R=2.pc @C=4.0pc 
{
(Y_\mu:X^\lambda)^\vartriangle \ar[d]_{\vartriangle_{Y_\mu:X^\lambda}} \ar[r]^{(\wedge_\mu^Y:\vee_\lambda^X)^\vartriangle} & (Y:X)^\vartriangle \ar@{-->}[d]^{h} \\
Y_\mu:X^\lambda \ar[r]_{g:x}& \C
}
\eeq
}\eit
\elm
\bpr
Здесь нужно провести двойную индукцию по $\lambda$ и $\mu$. При $\lambda=\mu=0$ равенство \eqref{(Y_mu:X_lambda)^vartriangle=(Y:X)^vartriangle-func} становится тривиальным
$$
(Y:X)^\vartriangle=(Y:X)^\vartriangle,
$$
а диаграмма \eqref{(Y_mu:X_lambda)^vartriangle=(Y:X)^vartriangle-func} превращается в диаграмму
$$
\xymatrix @R=2.pc @C=8.0pc 
{
(Y:X)^\vartriangle \ar[d]_{\vartriangle_{Y:X}} \ar[r]^{(\id_Y:\id_X)^\vartriangle=\id_{Y:X}^\vartriangle} & (Y:X)^\vartriangle \ar@{-->}[d]^{h} \\
Y:X \ar[r]_{g:x}& \C
}
$$
в которой нужно просто положить $h=g:x\ \circ\vartriangle_{X\cdot Y}$.

Предположим, что мы доказали утверждения $(i)^\circ$ и $(ii)^\circ$ для всех ординалов $\iota<\lambda$ и $\varkappa<\mu$, где $\lambda$ и $\mu$ --- какие-то два фиксированных ординала:
\bit{
\item[$(i)^\circ$] при псевдонасыщении пространства $Y_\varkappa:X^\iota$ и $Y:X$ становятся изоморфными:
\beq\label{(Y_varkappa:X_iota)^vartriangle=(Y:X)^vartriangle}
(Y_\varkappa:X^\iota)^\vartriangle=(Y:X)^\vartriangle,
\eeq

\item[$(ii)^\circ$] для любых $x\in X^\iota$ и $g\in (Y_\varkappa)^\star$ функционал
$$
    (g:x):(Y_\varkappa:X^\iota)\to\C\quad\Big|\quad (g:x)(\ph)=g(\ph(x)),\qquad \ph\in Y_\varkappa:X^\iota
$$
    непрерывен на пространстве $(Y:X)^\vartriangle$, то есть  существует единственный функционал $h\in ((Y:X)^\vartriangle)^\star$, замыкающий диаграмму
\beq\label{(Y_varkappa:X_iota)^vartriangle=(Y:X)^vartriangle-func}
\xymatrix @R=2.pc @C=4.0pc 
{
(Y_\varkappa:X^\iota)^\vartriangle \ar[d]_{\vartriangle_{Y_\varkappa:X^\iota}} \ar[r]^{(\wedge_\varkappa^Y:\vee_\iota^X)^\vartriangle} & (Y:X)^\vartriangle \ar@{-->}[d]^{h} \\
Y_\varkappa:X^\iota \ar[r]_{g:x}& \C
}
\eeq
}\eit

Покажем, что тогда (i) и (ii) выполняются для ординалов $\lambda$ и $\mu$.

1. Сначала зафиксируем произвольный ординал $\varkappa$ так, чтобы $\varkappa<\mu$ и покажем, что $(i)^\circ$ и $(ii)^\circ$ выполняются при подстановке $\iota=\lambda$. Здесь нужно рассмотреть два случая.

a) Предположим сначала, что $\lambda$ -- изолированный ординал, то есть
$$
\lambda=\iota+1
$$
для некоторого ординала $\iota<\lambda$. Покажем, что тогда для $\lambda$ (подставленного вместо $\iota$) выполняется условие $(ii)^\circ$. Пусть
$$
x\in X^\lambda=X^{\iota+1}=(X^\iota)^\vee
$$
и
$$
g\in (Y_\varkappa)^\star
$$
Тогда $x$ является пределом некоторой вполне ограниченной направленности $\{x_i;\ i\in I\}\subseteq X^\iota$:
\beq\label{(Y_mu:X_lambda)^vartriangle=(Y:X)^vartriangle-1}
x_i\overset{(X^\iota)^\vee}{\underset{i\to\infty}{\longrightarrow}}x
\eeq
Поскольку множество $\{x_i\}\subseteq X^\iota$ вполне ограничено, по следствию \ref{COR:(g,x)|->g:x-razdelno-nepreryvno}, множество
$$
\{g:x_i;\ i\in I\}
$$
вполне ограничено в $(Y_\varkappa:X^\iota)^\star$. С другой стороны, из теоремы \ref{TH:(g,x)|->g:x-razdelno-nepreryvno} следует, что $\{g:x_i;\ i\in I\}$ является направленностью Коши. Итак, мы получаем, что $\{g:x_i;\ i\in I\}$ --- вполне ограниченная направленность Коши в пространстве $(Y_\varkappa:X^\iota)^\star$, значит, она сходится в объемлющем псевдополном пространстве $(Y_\varkappa:X^\iota)^{\star\triangledown}$:
$$
(Y_\varkappa:X^\iota)^\star\subseteq (Y_\varkappa:X^\iota)^{\star\triangledown}=
\text{\cite[Theorem 3.14]{Ak03}}=(Y_\varkappa:X^\iota)^{\vartriangle\star}= \eqref{(Y_varkappa:X_iota)^vartriangle=(Y:X)^vartriangle}=
(Y:X)^{\vartriangle\star}
$$
То есть найдется функционал $h\in (Y:X)^{\vartriangle\star}$ такой что
\beq\label{(Y_mu:X_lambda)^vartriangle=(Y:X)^vartriangle-2}
g:x_i\overset{(Y:X)^{\vartriangle\star}}{\underset{i\to\infty}{\longrightarrow}}h
\eeq
С другой стороны, из цепочки равенств
$$
(X^\iota)^\vee=X^{\iota+1}=X^\lambda
$$
мы получаем, что соотношение \eqref{(Y_mu:X_lambda)^vartriangle=(Y:X)^vartriangle-1} эквивалентно соотношению
\beq\label{(Y_mu:X_lambda)^vartriangle=(Y:X)^vartriangle-3}
x_i\overset{X^\lambda}{\underset{i\to\infty}{\longrightarrow}} x
\eeq
из которого по теореме \ref{TH:(g,x)|->g:x-razdelno-nepreryvno} следует
\beq\label{(Y_mu:X_lambda)^vartriangle=(Y:X)^vartriangle-4}
g:x_i\overset{(Y_\varkappa:X^\lambda)^{\star}}{\underset{i\to\infty}{\longrightarrow}} g:x
\eeq
Вместе соотношения \eqref{(Y_mu:X_lambda)^vartriangle=(Y:X)^vartriangle-2} и \eqref{(Y_mu:X_lambda)^vartriangle=(Y:X)^vartriangle-4} означают, в частности, что на каждом операторе
$$
\ph\in Y_\varkappa:X^\lambda
$$
(по лемме \ref{LM:Y_mu:X^lambda=Y:X}, пространства $Y_\varkappa:X^\lambda$ и $Y:X$ совпадают как множества, различаясь только топологией) функционалы $g:x$ и $h$ совпадают:
$$
(g:x)(\ph)=h(\ph).
$$
Это доказывает утверждение $(ii)^\circ$ (для случая $\lambda=\iota+1$ подставленного вместо $\iota$).

Докажем теперь $(i)^\circ$ (опять для случая $\lambda=\iota+1$ подставленного вместо $\iota$). Отображение
$$
(Y_\varkappa:X^\lambda)^\vartriangle\to (Y:X)^\vartriangle
$$
в прямую сторону всегда непрерывно. Нужно доказать, что оно непрерывно в обратную сторону, то есть, что непрерывно отображение
$$
(Y:X)^\vartriangle\to (Y_\varkappa:X^\lambda)^\vartriangle.
$$
Поскольку пространство $(Y:X)^\vartriangle$ псевдонасыщено, нам достаточно доказать непрерывность отображения
$$
(Y:X)^\vartriangle\to (Y_\varkappa:X^\lambda).
$$
Рассмотрим базисную окрестность нуля в $Y_\varkappa:X^\lambda$, то есть множество $V:S$, где $V$ --- замкнутая выпуклая уравновешенная окрестность нуля в $Y_\varkappa$, а $S$ --- вполне ограниченное множенство в $X^\lambda$. По лемме \ref{LM:o-drobyah-v-Y_mu:X_lambda}, множество $V:S$ является емким в $Y:X$. Поэтому если мы докажем, что оно замкнуто в $(Y:X)^\vartriangle$, то это будет означать, что оно является окрестностью нуля в $(Y:X)^\vartriangle$.

Это становится очевидно, если представить $V:S$ как поляру системы функционалов вида $\{g:x;\ g\in V^\circ,\ x\in S\}$:
\beq\label{(Y_mu:X_lambda)^vartriangle=(Y:X)^vartriangle-5}
V:S=\{g:x;\ g\in V^\circ,\ x\in S\}^\circ
\eeq
Действительно,
\begin{multline}\label{(Y_mu:X_lambda)^vartriangle=(Y:X)^vartriangle-5-proof}
\ph\in V:S \quad\Leftrightarrow\quad \ph(S)\subseteq V
\quad\Leftrightarrow\quad \forall x\in S \ \ph(x)\in V
\quad\Leftrightarrow \\ \Leftrightarrow\quad
\forall g\in V^\circ\ \forall x\in S \ |\underbrace{g(\ph(x))}_{\scriptsize \begin{matrix}\|\\ (g:x)(\ph)\end{matrix}}|\le 1
\quad\Leftrightarrow\quad \ph\in \{g:x;\ g\in V^\circ,\ x\in S\}^\circ
\end{multline}
Теперь нужно заметить, что в формуле \eqref{(Y_mu:X_lambda)^vartriangle=(Y:X)^vartriangle-5}
$x\in S$, где $S$ --- вполне ограниченное множество в $X^\lambda$, а $g\in V^\circ$, где $V$ --- окрестность нуля в $Y_\varkappa$. Значит, $x\in X^\lambda$, а  $g\in (Y_\varkappa)^\star$. Отсюда по уже доказанному свойству $(ii)^\circ$ (для случая $\lambda=\iota+1$),  $g:x$ --- непрерывные функционалы на пространстве $(Y:X)^\vartriangle$. Мы получаем, что $V:S$ является полярой некоторой системы непрерывных функционалов на $(Y:X)^\vartriangle$, значит, $V:S$ --- замкнутое множество в $(Y:X)^\vartriangle$. Вдобавок, оно емкое, значит, оно является окрестностью нуля в (псевдонасыщенном) пространстве $(Y:X)^\vartriangle$. Мы доказали свойство $(i)^\circ$ (для случая $\lambda=\iota+1$ подставленного вместо $\iota$).

b) Теперь рассмотрим случай, когда $\lambda$ --- предельный ординал, то есть
$$
\lambda\ne \iota+1
$$
для всех $\iota<\lambda$. Здесь опять сначала нужно доказать $(ii)^\circ$ (для $\lambda$ подставленного вместо $\iota<\lambda$). Пусть $x\in X^\lambda$ и $g\in (Y_\varkappa)^\star$. Тогда
$$
x\in X^\lambda=\text{\cite[(1.19)]{Ak03}}=\bigcup_{\iota<\lambda} X^\iota,
$$
поэтому существует ординал $\iota<\lambda$ такой, что $x\in X^\iota$. Мы получаем, что $x\in X^\iota$ и $g\in (Y_\varkappa)^\star$, поэтому по предположению индукции $(ii)^{\circ}$, функционал $g:x$ должен быть непрывен на пространстве $(Y:X)^\vartriangle$. Это доказывает  $(ii)^\circ$ (для случая предельного $\lambda$ подставленного вместо $\iota<\lambda$). Теперь перейдем к $(i)^\circ$. Здесь повторяются рассуждения пункта a): нам нужно доказать непрерывность отображения
$$
(Y:X)^\vartriangle\to (Y_\varkappa:X^\lambda).
$$
Мы берем базисную окрестность нуля в $Y_\varkappa:X^\lambda$, то есть множество $V:S$, где $S$ --- вполне ограниченное множество в $X^\lambda$, а $V$ --- замкнутая выпуклая уравновешенная окрестность нуля в $Y_\varkappa$, затем доказываем равенство
\eqref{(Y_mu:X_lambda)^vartriangle=(Y:X)^vartriangle-5}, затем замечаем, что в нем $x\in X^\lambda$, а  $g\in (Y_\varkappa)^\star$. После этого по уже доказанному свойству $(ii)^\circ$ (для случая предельного $\lambda$ подставленного вместо $\iota<\lambda$), мы получаем, что $g:x$ --- непрерывные функционалы на пространстве $(Y:X)^\vartriangle$. Значит, $V:S$ --- замкнутое множество в $(Y:X)^\vartriangle$. Вдобавок, оно емкое, по лемме \ref{LM:o-drobyah-v-Y_mu:X_lambda}, значит, оно является окрестностью нуля в (псевдонасыщенном) пространстве $(Y:X)^\vartriangle$. Мы доказали свойство $(i)^\circ$ (для случая предельного $\lambda$ подставленного вместо $\iota<\lambda$).

2. Мы показали, что при фиксированном $\varkappa<\mu$ условия $(i)^\circ$ и $(ii)^\circ$ остаются справедливыми, если мы в них перейдем от $\iota<\lambda$ к $\iota=\lambda$. Теперь зафиксируем ординал $\iota$, положив $\iota=\lambda$, и покажем, что $(i)^\circ$ и $(ii)^\circ$ (с подставленным $\iota=\lambda$) остаются справедливыми, если в них мы  перейдем от $\varkappa<\mu$ к $\varkappa=\mu$. Здесь тоже нужно рассмотреть два случая.

a) Предположим сначала, что $\mu$ -- изолированный ординал, то есть
$$
\mu=\varkappa+1
$$
для некоторого ординала $\varkappa<\mu$. Покажем, что тогда для $\mu$ (подставленного вместо $\varkappa$) выполняется условие $(ii)^\circ$. Пусть
$$
g\in (Y_\mu)^\star=(Y_{\varkappa+1})^\star=((Y_\varkappa)^\wedge)^\star
=\text{\cite[Theorem 3.10]{Ak03}}=((Y_\varkappa)^\star)^\vee
$$
и
$$
x\in X^\lambda.
$$
Тогда $g$ является пределом некоторой вполне ограниченной направленности $\{g_i;\ i\in I\}\subseteq (Y_\varkappa)^\star$:
\beq\label{(Y_mu:X_lambda)^vartriangle=(Y:X)^vartriangle-1*}
g_i\overset{((Y_\varkappa)^\star)^\vee}{\underset{i\to\infty}{\longrightarrow}}g
\eeq
Поскольку множество $\{g_i\}\subseteq (Y_\varkappa)^\star$ вполне ограничено, по следствию \ref{COR:(g,x)|->g:x-razdelno-nepreryvno}, множество
$$
\{g_i:x;\ i\in I\}
$$
вполне ограничено в $(Y_\varkappa:X^\lambda)^\star$. С другой стороны, из теоремы \ref{TH:(g,x)|->g:x-razdelno-nepreryvno} следует, что $\{g_i:x;\ i\in I\}$ является направленностью Коши. Итак, мы получаем, что $\{g_i:x;\ i\in I\}$ --- вполне ограниченная направленность Коши в пространстве $(Y_\varkappa:X^\lambda)^\star$, значит, она сходится в объемлющем псевдополном пространстве $(Y_\varkappa:X^\lambda)^{\star\triangledown}$:
$$
(Y_\varkappa:X^\lambda)^\star\subseteq (Y_\varkappa:X^\lambda)^{\star\triangledown}=
\text{\cite[Theorem 3.14]{Ak03}}=(Y_\varkappa:X^\lambda)^{\vartriangle\star}= \eqref{(Y_varkappa:X_iota)^vartriangle=(Y:X)^vartriangle}=
(Y:X)^{\vartriangle\star}
$$
То есть найдется функционал $h\in (Y:X)^{\vartriangle\star}$ такой что
\beq\label{(Y_mu:X_lambda)^vartriangle=(Y:X)^vartriangle-2*}
g_i:x\overset{(Y:X)^{\vartriangle\star}}{\underset{i\to\infty}{\longrightarrow}}h
\eeq
С другой стороны, из цепочки равенств
$$
((Y_\varkappa)^\star)^\vee= \text{\cite[Theorem 3.10]{Ak03}}= ((Y_\varkappa)^\wedge)^\star=(Y_{\varkappa+1})^\star=(Y_\mu)^\star
$$
мы получаем, что соотношение \eqref{(Y_mu:X_lambda)^vartriangle=(Y:X)^vartriangle-1*} эквивалентно соотношению
\beq\label{(Y_mu:X_lambda)^vartriangle=(Y:X)^vartriangle-3*}
g_i\overset{(Y_\mu)^\star}{\underset{i\to\infty}{\longrightarrow}} g
\eeq
из которого по теореме \ref{TH:(g,x)|->g:x-razdelno-nepreryvno} следует
\beq\label{(Y_mu:X_lambda)^vartriangle=(Y:X)^vartriangle-4*}
g_i:x\overset{(Y_\mu:X^\lambda)^{\star}}{\underset{i\to\infty}{\longrightarrow}} g:x
\eeq
Вместе соотношения \eqref{(Y_mu:X_lambda)^vartriangle=(Y:X)^vartriangle-2*} и \eqref{(Y_mu:X_lambda)^vartriangle=(Y:X)^vartriangle-4*} означают, в частности, что на каждом операторе
$$
\ph\in Y_\mu:X^\lambda
$$
(по лемме \ref{LM:Y_mu:X^lambda=Y:X}, пространства $Y_\mu:X^\lambda$ и $Y:X$ совпадают как множества, различаясь только топологией) функционалы $g:x$ и $h$ совпадают:
$$
(g:x)(\ph)=h(\ph).
$$
Это доказывает утверждение $(ii)^\circ$ (для случая $\mu=\varkappa+1$ подставленного вместо $\varkappa$).

Докажем теперь $(i)^\circ$ (опять для случая $\mu=\varkappa+1$ подставленного вместо $\varkappa$). Отображение
$$
(Y_\mu:X^\lambda)^\vartriangle\to (Y:X)^\vartriangle
$$
в прямую сторону всегда непрерывно. Нужно доказать, что оно непрерывно в обратную сторону, то есть, что непрерывно отображение
$$
(Y:X)^\vartriangle\to (Y_\mu:X^\lambda)^\vartriangle.
$$
Поскольку пространство $(Y:X)^\vartriangle$ псевдонасыщено, нам достаточно доказать непрерывность отображения
$$
(Y:X)^\vartriangle\to (Y_\mu:X^\lambda).
$$
Рассмотрим базисную окрестность нуля в $Y_\mu:X^\lambda$, то есть множество $V:S$, где $V$ --- замкнутая выпуклая уравновешенная окрестность нуля в $Y_\mu$, а $S$ --- вполне ограниченное множество в $X^\lambda$. По лемме \ref{LM:o-drobyah-v-Y_mu:X_lambda}, множество $V:S$ является емким в $Y:X$. Поэтому если мы докажем, что оно замкнуто в $(Y:X)^\vartriangle$, то это будет означать, что оно является окрестностью нуля в $(Y:X)^\vartriangle$.

Это становится очевидно, если представить $V:S$ как поляру системы функционалов вида $\{g:x;\ g\in V^\circ,\ x\in S\}$:
\beq\label{(Y_mu:X_lambda)^vartriangle=(Y:X)^vartriangle-5*}
V:S=\{g:x;\ g\in V^\circ,\ x\in S\}^\circ
\eeq
что доказывается той же цепочкой \eqref{(Y_mu:X_lambda)^vartriangle=(Y:X)^vartriangle-5-proof}, что и раньше.
Далее мы замечаем, что в формуле \eqref{(Y_mu:X_lambda)^vartriangle=(Y:X)^vartriangle-5*}
$x\in S$, где $S$ --- вполне ограниченное множество в $X^\lambda$, а $g\in V^\circ$, где $V$ --- окрестность нуля в $Y_\mu$. Значит, $x\in X^\lambda$, а  $g\in (Y_\mu)^\star$. Отсюда по уже доказанному свойству $(ii)^\circ$ (для случая $\mu=\varkappa+1$),  $g:x$ --- непрерывные функционалы на пространстве $(Y:X)^\vartriangle$. Мы получаем, что $V:S$ является полярой некоторой системы непрерывных функционалов на $(Y:X)^\vartriangle$, значит, $V:S$ --- замкнутое множество в $(Y:X)^\vartriangle$. Вдобавок, оно емкое, значит, оно является окрестностью нуля в (псевдонасыщенном) пространстве $(Y:X)^\vartriangle$. Мы доказали свойство $(i)^\circ$ (для случая $\mu=\varkappa+1$ подставленного вместо $\varkappa$).

b) Теперь рассмотрим случай, когда $\mu$ --- предельный ординал, то есть
$$
\mu\ne \varkappa+1
$$
для всех $\varkappa<\mu$. Здесь опять сначала нужно доказать $(ii)^\circ$ (с подставленным $\mu$ вместо $\varkappa<\mu$). Пусть $x\in X^\lambda$ и $g\in (Y_\mu)^\star$. Тогда
$$
g\in (Y_\mu)^\star=\text{\cite[Theorem 3.12]{Ak03}} =(Y^\star)^\mu=\text{\cite[(1.19)]{Ak03}}=\bigcup_{\varkappa<\mu} (Y^\star)^\varkappa=\text{\cite[Theorem 3.12]{Ak03}} =\bigcup_{\varkappa<\mu}(Y_\varkappa)^\star,
$$
поэтому существует ординал $\varkappa<\mu$ такой, что $g\in (Y^\star)^\varkappa=(Y_\varkappa)^\star$. Мы получаем, что $x\in X^\iota$ и $g\in (Y_\varkappa)^\star$, поэтому по предположению индукции $(ii)^{\circ}$, функционал $g:x$ должен быть непрывен на пространстве $(Y:X)^\vartriangle$. Это доказывает  $(ii)^\circ$ (для случая предельного $\mu$ подставленного вместо $\varkappa<\mu$). Теперь перейдем к $(i)^\circ$. Здесь повторяются рассуждения пункта a): нам нужно доказать непрерывность отображения
$$
(Y:X)^\vartriangle\to (Y_\mu:X^\lambda).
$$
Мы берем базисную окрестность нуля в $Y_\mu:X^\lambda$, то есть множество $V:S$, где $S$ --- вполне ограниченное множество в $X^\lambda$, а $V$ --- замкнутая выпуклая уравновешенная окрестность нуля в $Y_\mu$, затем доказываем равенство
\eqref{(Y_mu:X_lambda)^vartriangle=(Y:X)^vartriangle-5*}, затем замечаем, что в нем $x\in X^\lambda$, а  $g\in (Y_\mu)^\star$. После этого по уже доказанному свойству $(ii)^\circ$ (для случая предельного $\mu$ подставленного вместо $\varkappa<\mu$), мы получаем, что $g:x$ --- непрерывные функционалы на пространстве $(Y:X)^\vartriangle$. Значит, $V:S$ --- замкнутое множество в $(Y:X)^\vartriangle$. Вдобавок, оно емкое, по лемме \ref{LM:o-drobyah-v-Y_mu:X_lambda}, значит, оно является окрестностью нуля в (псевдонасыщенном) пространстве $(Y:X)^\vartriangle$. Мы доказали свойство $(i)^\circ$ (для случая предельного $\mu$ подставленного вместо $\varkappa<\mu$).
\epr

Для произвольных локально выпуклых пространств $X$ и $Y$ отображения
$$
\triangledown_X: X\to X^\triangledown,\qquad \vartriangle_Y: Y^\vartriangle\to Y
$$
определяют отображение
\beq\label{(vartriangle_Y:triangledown_X)^vartriangle}
(\vartriangle_Y:\triangledown_X): (Y^\vartriangle:X^\triangledown) \to (Y:X).
\eeq
(которое, понятное дело, просто совпадает с отображением \eqref{wedge_mu^Y:vee_lambda^X-1} при подходящих $\lambda$ и $\mu$).

\btm\label{TH:(vartriangle_Y:triangledown_X)^vartriangle-isomorphism}
Пусть $X$ псевдонасыщено, а $Y$ псевдополно. Тогда псевдонасыщение отображения \eqref{(vartriangle_Y:triangledown_X)^vartriangle} является изоморфизмом локально выпуклых (и стереотипных) пространств:
\beq\label{(vartriangle_Y:triangledown_X)^vartriangle-isomorphism}
Y^\vartriangle\oslash X^\triangledown\cong (Y^\vartriangle:X^\triangledown)^\vartriangle\cong  (Y:X)^\vartriangle
\eeq
\etm
\bpr
Здесь второе равенство получается из \eqref{(Y_mu:X_lambda)^vartriangle=(Y:X)^vartriangle}, если подобрать ординалы $\lambda$ и $\mu$ так, чтобы $X^\lambda=X^\triangledown$ и $Y_\mu=Y^\vartriangle$.
А первое --- просто определение \eqref{DEF:oslash} тензорной дроби $\oslash$ (при этом пространства $X^\triangledown$ и $Y^\vartriangle$ стереотипны в силу \cite[Proposition 3.17]{Ak03} и \cite[Proposition 3.16]{Ak03}).
\epr

\section{Псевдонасыщение первичного тензорного произведения $X\cdot Y$}

\subsection{Первичное тензорное произведение $X\cdot Y$.}

Пусть $X$ и $Y$ --- локально выпуклые пространства. Условимся {\it первичным тензорным произведением} $X\cdot Y$ называть локально выпуклое пространство, состоящее из операторов $\ph:X^\star\to Y$, и наделенное {\it топологией равномерной сходимости на полярах окрестностей нуля $U\subseteq X$}:
\beq\label{shodimost-v-X-cdot-Y}
\ph_i\overset{X\cdot Y}{\underset{i\to\infty}{\longrightarrow}}\ph\quad\Leftrightarrow\quad
\forall U\in{\mathcal U}(X)\quad
\ph_i(f)\overset{Y}{\underset{f\in U^\circ}{\underset{i\to\infty}{\rightrightarrows}}}\ph(f)
\eeq
Эту топологию удобно обозначить какой-нибудь буквой, например, $\xi$, тогда пространство $X\cdot Y$ можно будет представить формулой
\beq\label{DEF:X-cdot-Y}
X\cdot Y=Y\underset{\xi}{:} X^\star
\eeq
(индекс $\xi$ обозначает сходимость в топологии $\xi$). Очевидно, что существует биективный оператор
$$
Y:X\to Y\underset{\xi}{:} X^\star=X\cdot Y
$$
(который, однако, не является изоморфизмом, если $X$ не псевдонасыщено).
Если $A\subseteq X$ и $B\subseteq Y$, то символ $A\cdot B$ обозначает множество
$$
A\cdot B=B:(A^\circ)
$$
(где двоеточие : было определено в \eqref{DEF:B:A}).

Если $\alpha:X'\to X$, $\beta:Y'\to Y$ --- операторы, то диаграмма
\beq\label{alpha-cdot-beta}
\xymatrix @R=2.pc @C=2.0pc 
{
X^\star \ar[d]_{(\alpha\cdot\beta)(\psi)=\beta\circ\psi\circ\alpha^\star} \ar[r]^{\alpha^\star} & (X')^\star \ar[d]^{\psi} \\
Y & Y' \ar[l]_{\beta}
}
\eeq
определяет отображение
$$
\ph\cdot\chi: X'\cdot Y'\to X\cdot Y
$$

\bprop\label{alpha-cdot-beta-nepreryvnoe-otobr}
Для любых операторов $\alpha:X'\to X$ и $\beta:Y'\to Y$ отображение
$$
\ph\cdot\chi: X'\cdot Y'\to X\cdot Y
$$
непрерывно.
\eprop
\bpr Из непрерывности $\alpha$ следует, что сопряженное отображение $\alpha^\star:X^\star\to (X')^\star$ переводит поляру любой окрестности нуля $U\subseteq X$ в поляру некоторой окрестности нуля $U'$, а именно, $U'=\alpha^{-1}(U)$:
\begin{multline*}
f\in U^\circ\quad\Rightarrow\quad \forall x\in U\ \abs{f(x)}\le 1
\quad\Rightarrow\\ \Rightarrow\quad \forall x'\in U'=\alpha^{-1}(U)\ \abs{\alpha^\star(f)(x')}=|f(\underbrace{\alpha(x')}_{\scriptsize \begin{matrix}\text{\rotatebox{90}{$\owns$}}\\ U\end{matrix}})|\le 1
\quad\Rightarrow\quad \alpha^\star(f)\in (U')^\circ
\end{multline*}
Поэтому из сходимости $\psi_i\to \psi$ в $X'\cdot Y'$
$$
\forall U'\in{\mathcal U}(X')\quad
\psi_i(g)\overset{Y'}{\underset{g\in (U')^\circ}{\underset{i\to\infty}{\rightrightarrows}}}\psi(g)
$$
следует сходимость $(\alpha\cdot\beta)(\psi_i)\to (\alpha\cdot\beta)(\psi)$ в $X\cdot Y$:
$$
\forall U\in{\mathcal U}(X)\quad
(\alpha\cdot\beta)(\psi_i)(f)=\beta\Big(\psi_i\big(\alpha^\star(f)\big)\Big)
\overset{Y}{\underset{f\in U^\circ}{\underset{i\to\infty}{\rightrightarrows}}}
\beta\Big(\psi\big(\alpha^\star(f)\big)\Big)=
(\alpha\cdot\beta)(\psi)(f)
$$
\epr

Пусть далее $X$ и $Y$ --- локально выпуклые пространства, причем $X$ псевдополно. Для всякого оператора $\ph:X^\star\to Y$ мы можем рассмотреть двойственный оператор $\ph:Y^\star\to X^{\star\star}$, а после этого сформировать композицию с оператором $i_X^{-1}:X^{\star\star}\to X$ (который в силу \cite[Corollary 2.13]{Ak03} существует и непрерывен, потому что $X$ псевдополно). Обозначим эту композицию символом $\omega_{X,Y}(\ph)$:
\beq\label{DEF:omega(ph)}
\xymatrix @R=2.pc @C=6.0pc 
{
Y^\star \ar@/_4ex/[rr]_{\omega_{X,Y}(\ph)} \ar[r]^{\ph^\star} & X^{\star\star} \ar[r]^{i_X^{-1}} & X
}
\eeq
Мы получаем отображение
\beq\label{omega:XY->YX}
\omega_{X,Y}:X\cdot Y=Y\underset{\xi}{:}X^\star\to X\underset{\xi}{:}Y^\star=Y\cdot X\quad\Big|\quad \omega_{X,Y}(\ph)=i_X^{-1}\circ\ph^\star
\eeq

\btm\label{TH:X-cdot-Y-cong-Y-cdot-X}
Если локально выпуклые пространства $X$ и $Y$ псевдополны, то отображение \eqref{omega:XY->YX} устанавливает изоморфизм локально выпуклых пространств
\beq\label{X-cdot-Y-cong-Y-cdot-X}
X\cdot Y\cong Y\cdot X
\eeq
\etm
\bpr

1. Сначала заметим такое тождество:
\beq\label{f(i_X^-1(ph^star(g)))=g((ph(f)))}
f(\omega_{X,Y}(\ph)(g))=g(\ph(f)),\qquad f\in X^\star, \ g\in Y^\star,\ \ph\in Y:(X^\star).
\eeq
Для доказательства нужно обозначить
$$
x=\omega_{X,Y}(\ph)(g)=i_X^{-1}(\ph^\star(g))=i_X^{-1}(g\circ\ph)\in X,
$$
тогда мы получим
$$
f(\underbrace{\omega_{X,Y}(\ph)(g)}_{\scriptsize\begin{matrix}\|\\ x\end{matrix}})=f(x)=i_X(x)(f)=
i_X(i_X^{-1}(g\circ\ph))(f)=(g\circ\ph)(f)=g(\ph(f))
$$

2. После того, как доказано \eqref{f(i_X^-1(ph^star(g)))=g((ph(f)))}, убедимся, что отображение $\omega_{X,Y}$ биективно. Это следует из тождества
\beq
\omega_{Y,X}(\omega_{X,Y}(\ph))=\ph,\qquad \ph\in X\cdot Y=Y\underset{\xi}{:}X^\star.
\eeq
Действительно, для любых $f\in X^\star$ и $g\in Y^\star$ мы получим
$$
g(\omega_{Y,X}(\omega_{X,Y}(\ph))(f))=\eqref{f(i_X^-1(ph^star(g)))=g((ph(f)))}=
f(\omega_{X,Y}(\ph)(g))=\eqref{f(i_X^-1(ph^star(g)))=g((ph(f)))}=g(\ph(f))
$$

3. Наконец, докажем непрерывность отображение $\omega_{X,Y}$  в обе стороны. Действительно, если $U:(V^\circ)$ --- базисная окрестность нуля в $X\underset{\xi}{:}Y^\star=Y\cdot X$ (где $U\subseteq X$ и $V\subseteq Y$ --- замкнутые выпуклые уравновешенные окрестности нуля), то $V:(U^\circ)$ --- базисная окрестность нуля в $Y\underset{\xi}{:}X^\star=X\cdot Y$, и при отображениях $\omega_{X,Y}$ и $\omega_{X,Y}^{-1}$ одна окрестность переходит в другую:
\begin{multline*}
\omega_{X,Y}(\ph)\in U:(V^\circ)\quad\Leftrightarrow\quad \omega_{X,Y}(\ph)(V^\circ)\subseteq U
\quad\Leftrightarrow\quad \forall f\in U^\circ\ \forall g\in V^\circ\quad |\underbrace{f(\omega_{X,Y}(\ph)(g))}_{\scriptsize \begin{matrix}\|\put(3,0){\eqref{f(i_X^-1(ph^star(g)))=g((ph(f)))}}
\\ g(\ph(f))\end{matrix}}|\le 1
\quad\Leftrightarrow
\\ \Leftrightarrow \quad
\forall f\in U^\circ\ \forall g\in V^\circ\quad |g(\ph(f))|\le 1
\quad\Leftrightarrow\quad \ph(U^\circ)\subseteq V \quad\Leftrightarrow\quad \ph\in V:(U^\circ)
\end{multline*}
\epr

Вернемся к диаграмме \eqref{alpha-cdot-beta}. В частном случае, когда $\alpha$ и $\beta$ представляют собой функционалы
$$
f:X\to\C,\qquad g:Y\to\C,
$$
отображение $\alpha\cdot\beta$ также естественно понимать как функционал
$$
f\cdot g:X\cdot Y\to\C,
$$
действующий по формуле
\beq\label{DEF:f-cdot-g}
(f\cdot g)(\ph)=g(\ph(f)),\qquad \ph:X^\star\to Y.
\eeq

Если $F\subseteq X^\star$ и $G\subseteq Y^\star$ -- два множества функционалов, то
символом $F\cdot G$ мы обозначаем множество
$$
F\cdot G=\{f\cdot g;\ f\in F,\ g\in G\}.
$$
В частности, если $f\in X^\star$ и $G\subseteq Y^\star$, то
$$
f\cdot G=\{f\cdot g;\ g\in G\}.
$$
а если $F\subseteq X^\star$ и $g\in Y^\star$, то
$$
F\cdot g=\{f\cdot g;\ f\in F\}.
$$

\btm\label{TH:(f,g)|->f-cdot-g-razdelno-nepreryvno}
Для псевдополных локально выпуклых пространств $X$ и $Y$ отображение
$$
(f,g)\in X^\star\times Y^\star\mapsto f\cdot g\in (X\cdot Y)^\star
$$
раздельно непрерывно.
\etm
\bpr
Пусть $\varPhi\subseteq X\cdot Y=Y\underset{\xi}{:}X^\star$ --- вполне ограниченное множество и $f\in X^\star$. По определению топологии $\xi$, $\varPhi$ должно быть равностепенно непрерывно и равномерно вполне ограничено на поляре $U^\circ$ любой окрестности нуля $U\subseteq X$ (здесь применяется вариант теоремы 5.1 из \cite{Ak03}). В частности, если взять окрестность $U=\{f\}^\circ$, то множество операторов $\varPhi$ будет переводить его поляру
$$
U^\circ=\{f\}^{\circ\circ}=\{\lambda\in\C:\ |\lambda|\le 1\}\cdot f.
$$
во вполне ограниченное множество в пространстве $Y$:
$$
\varPhi(U^\circ)=\{\lambda\in\C:\ |\lambda|\le 1\}\cdot \varPhi(f)\subseteq Y.
$$
Как следствие, поляра этого множества
$$
V=\{\lambda\in\C:\ |\lambda|\le 1\}\cdot \varPhi(f)^\circ\subseteq Y^\star
$$
должна быть окрестностью нуля в $Y^\star$. Теперь мы получаем, что для любого $g\in V$
$$
\sup_{\ph\in\varPhi}\abs{(f\cdot g)(\ph)}=\sup_{\ph\in\varPhi}\abs{g(\ph(f))}\le 1
$$
то есть $f\cdot g\in\varPhi^\circ$. Это верно для любого $g\in V$, поэтому
\beq\label{(f,g)|->f-cdot-g-razdelno-nepreryvno-1}
f\cdot V\subseteq \varPhi^\circ.
\eeq

Иными словами, для всякой базисной окрестности нуля $\varPhi^\circ$ в $(X\cdot Y)^\star$ (где $\varPhi$ --- вполне ограниченное множество в $X\cdot Y$) и любого вектора $f\in X^\star$ найдется окрестность нуля $V\subseteq Y^\star$ такая, что выполняется \eqref{(f,g)|->f-cdot-g-razdelno-nepreryvno-1}. В силу \eqref{X-cdot-Y-cong-Y-cdot-X}, то же самое будет верно, если поменять местами $X$ и $Y$.
\epr

\bcor\label{LM:F-cdot-g-subseteq-(X-cdot-Y)^star}
Если $F\subseteq X^\star$ --- вполне ограниченное множество и $g\in Y^\star$, то  $F\cdot g\subseteq (X\cdot Y)^\star$ --- вполне ограниченное множество. Двойственным образом, если $f\in X^\star$ и $G\subseteq Y^\star$ --- вполне ограниченное множество, то $f\cdot G\subseteq (X\cdot Y)^\star$ --- вполне ограниченное множество.
\ecor

\bex\label{PROP:X-cdot-Y-cong-X-widetilde-otimes_e-Y}
Пусть $X$ и $Y$ --- полные локально выпуклые пространства, причем $Y$ обладает свойством (классической) аппроксимации. Тогда первичное тензорное произведение $X\cdot Y$ изоморфно инъективному тензорному произведению
\beq\label{X-cdot-Y-cong-X-widetilde-otimes_e-Y}
X\cdot Y\cong X\widetilde{\otimes}_{\e} Y.
\eeq
Если дополнительно $Y$ ядерно, то $X\cdot Y$ изоморфно также и проективному тензорному произведению
\beq\label{X-cdot-Y-cong-X-widetilde-otimes_pi-Y}
X\cdot Y\cong X\widetilde{\otimes}_{\e} Y\cong X\widetilde{\otimes}_{\pi} Y.
\eeq
\eex
\bpr
1. Рассмотрим вложение $X\otimes Y\subseteq X\cdot Y=Y\underset{\xi}{:} X^\star$, задаваемое формулой
$$
(x\otimes y)(f)=f(x)\cdot y,\qquad x\in X,\ y\in Y, \ f\in X^\star.
$$
Топология $\xi$ на $X\otimes Y$, индуцированная из $Y\underset{\xi}{:} X^\star$, совпадает с инъективной топологией $\e$, потому что $\xi$ и $\e$ порождаются одной и той же системой полунорм:
$$
\abs{z}_{U,V}=\sup_{f\in U^\circ,\ g\in V^\circ}\abs{g(z(f))}=
\sup_{f\in U^\circ,\ g\in V^\circ}\abs{g\bigg(\Big(\sum_{i=1}^n x_i\otimes y_i\Big)(f)\bigg)}=
\sup_{f\in U^\circ,\ g\in V^\circ}\abs{\sum_{i=1}^n f(x_i)\cdot g(y_i)}
$$
для $z=\sum_{i=1}^n x_i\otimes y_i\in X\otimes Y$.

С другой стороны, пространство $Y\underset{\xi}{:} X^\star$ полно, потому что $Y$ полно, а $X^\star$ насыщено (как сопряженное к полному $X$, в силу \cite[Theorem 2.17]{Ak03}). Отсюда следует, что инъективное тензорное произведение $X\widetilde{\otimes}_{\e} Y$ (то есть пополнение $X\otimes Y$ относительно топологии $\e$) представляет собой просто замыкание подпространства $X\otimes Y$ в пространстве $X\cdot Y=Y\underset{\xi}{:} X^\star$:
$$
X\widetilde{\otimes}_{\e} Y\cong\overline{X\otimes Y}.
$$
Чтобы доказать \eqref{X-cdot-Y-cong-X-widetilde-otimes_e-Y} нам нужно только убедиться, что $X\otimes Y$ плотно в $X\cdot Y=Y\underset{\xi}{:} X^\star$:
$$
\overline{X\otimes Y}=X\cdot Y.
$$
Это следует из того, что $Y$ обладает аппроксимацией: если $\psi_i$ --- направленность конечномерных операторов, сходящаяся к тождественному равномерно на компактах в $Y$,
$$
\psi_i\underset{i\to\infty}{\longrightarrow}\id_Y
$$
то для всякого $\ph\in X\cdot Y=Y\underset{\xi}{:} X^\star$ мы получим
$$
\psi_i\circ\ph\overset{X\cdot Y}{\underset{i\to\infty}{\longrightarrow}}\id_Y\circ\ph=\ph.
$$

2. Если $Y$ вдобавок ядерно, тогда второе равенство в \eqref{X-cdot-Y-cong-X-widetilde-otimes_pi-Y} получается автоматически из характеризации ядерных пространств \cite[7.3.3]{Pietsch}:
$$
X\widetilde{\otimes}_{\e} Y\cong X\widetilde{\otimes}_{\pi} Y.
$$
\epr

\subsection{Псевдонасыщение $X\cdot Y$.}

Пусть $\{X_\lambda;\ \lambda\in\bf{Ord}\}$ и $\{Y_\mu;\ \mu\in\bf{Ord}\}$ --- проективные ряды пространств $X$ и $Y$ \cite[1.4(b)]{Ak03}, и
$$
\wedge_\lambda^X: X_\lambda\to X,\qquad \wedge_\mu^Y: Y_\mu\to Y
$$
--- соответствующие естественные отображения. Рассмотрев аналог диаграммы \eqref{alpha-cdot-beta}
\beq\label{wedge_lambda^X-cdot-wedge_mu^Y}
\xymatrix @R=2.pc @C=2.0pc 
{
X^\star \ar[d]_{(\wedge_\lambda^X\cdot\wedge_\mu^Y)(\psi)= \wedge_\mu^Y\circ\psi\circ(\wedge_\lambda^X)^\star} \ar[r]^{(\wedge_\lambda^X)^\star} & (X_\lambda)^\star \ar[d]^{\psi} \\
Y & Y_\mu \ar[l]_{\wedge_\mu^Y}
}
\eeq
мы получим отображение
\beq\label{wedge_lambda^X-cdot-wedge_mu^Y-1}
\wedge_\lambda^X\cdot\wedge_\mu^Y: X_\lambda\cdot Y_\mu\to X\cdot Y,
\eeq
непрерывное по предложению \ref{alpha-cdot-beta-nepreryvnoe-otobr}.

\blm\label{LM:o-drobyah-*}
Если $X$ и $Y$ псевдополны, то отображение \eqref{wedge_lambda^X-cdot-wedge_mu^Y-1} является биекцией, поэтому пространства $X_\lambda\cdot Y_\mu$ и $X\cdot Y$ совпадают как множества (различаясь только топологией):
\beq\label{X_lambda-cdot-Y_mu=X-cdot-Y}
X_\lambda\cdot Y_\mu=X\cdot Y
\eeq
\elm
\bpr
Отображение $\wedge_\lambda^X: X_\lambda\to X$ является инъекцией, значит сопряженное отображение $(\wedge_\lambda^X)^\star: (X_\lambda)^\star\gets X^\star$ является эпиморфизмом. С другой стороны, отображение $\wedge_\mu^Y: Y_\mu\to Y$ является инъекцией. Отсюда следует, что отображение
$$
\psi\mapsto \wedge_\mu^Y\circ\psi\circ (\wedge_\lambda^X)^\star
$$
является инъекцией. Нам нужно проверить, что оно является сюръекцией. Действительно, пусть $\ph:X^\star\to Y$ --- произвольный оператор. Мы последовательно построим два оператора $\ph_\mu$ и $\psi$, замыкающие диаграмму:
$$
\xymatrix @R=3.pc @C=4.0pc 
{
X^\star \ar[d]_{\ph}\ar@{-->}[dr]_{\ph_\mu} \ar[r]^{(\wedge_\lambda^X)^\star} & (X_\lambda)^\star \ar@{-->}[d]^{\psi} \\
Y & Y_\mu \ar[l]^{\wedge_\mu^Y}
}
$$
Поскольку $X$ псевдополно, пространство $X^\star$ должно быть псевдонасыщено. Отсюда следует, что $\ph$ однозначно продолжается до оператора $\ph^\vartriangle:X^\star\to Y^\vartriangle$. Если теперь рассмотреть его композицию с естественным отображением $Y^\vartriangle\to Y_\mu$, то мы получим оператор, продолжающий $\ph$ до оператора $\ph_\mu:X^\star\to Y_\mu$.

Далее, поскольку $Y$ псевдополно, $Y_\mu$ тоже псевдополно. Поэтому оператор  $\ph_\mu:X^\star\to Y_\mu$ продолжается до некоторого оператора $\chi:(X^\star)^\triangledown\to Y_\mu$. Рассмотрев его композицию с естественным вложением $(X^\star)^\lambda\to (X^\star)^\triangledown$, мы получим оператор $\psi:(X^\star)^\lambda\to Y_\mu$. В силу равенства \cite[(3.11)]{Ak03}
$$
(X^\star)^\lambda=(X_\lambda)^\star,
$$
оператор $\psi$ можно считать оператором $\psi:(X_\lambda)^\star\to Y_\mu$.
\epr

\blm\label{LM:o-sohranenii-polnoi-ogranichennosti}
Пусть $X$ и $Y$ псевдополны. Если $\varPhi\in X\cdot Y$ --- вполне ограниченное множество, то для любых $\lambda,\mu\in\bf{Ord}$ его представление в пространстве $X_\lambda\cdot Y_\mu$ (биекцией \eqref{X_lambda-cdot-Y_mu=X-cdot-Y}) --- также вполне ограниченное множество.
\elm
\bpr
1. Рассмотрим сначала случай, когда $\lambda=0$. Если $\varPhi\subseteq Y\underset{\xi}{:}X^\star$ --- вполне ограниченное множество, то это значит, что  $\varPhi$ равностепенно непрерывно и равномерно вполне ограничено на поляре $U^\circ$ каждой окрестности нуля $U\subseteq X$ (мы используем здесь вариант \cite[Theorem 5.1]{Ak03}). Как следствие, образ $\varPhi(U^\circ)$ будет вполне ограниченным множеством в $Y$. Пространство $Y_\mu$ представляет собой усиление топологии на $Y$, при котором класс вполне ограниченных множеств и топология на вполне ограниченных множествах не меняются. Поэтому $\varPhi(U^\circ)$ будет вполне ограниченным множеством также и в пространстве $Y_\mu$, причем с той же топологией, что и индуцированная из $Y$. Отсюда можно заключить, что множество $\varPhi$, рассматриваемое как множество операторов из $X^\star$ в $Y_\mu$, равностепенно непрерывно и равномерно вполне ограничено на $U^\circ$. И это верно для всякой окрестности нуля $U\subseteq X$. Значит, $\varPhi$ вполне ограничено в $Y_\mu\underset{\xi}{:}X^\star=X\cdot Y_\mu$.

2. Итак, мы поняли, что если $\varPhi\subseteq X\cdot Y$ вполне ограничено в $X\cdot Y$, то оно вполне ограничено и в $X\cdot Y_\mu$ при любом $\mu\in\bf{Ord}$. По теореме \ref{TH:X-cdot-Y-cong-Y-cdot-X} отсюда следует, что $\varPhi$ вполне ограничено в пространстве $Y_\mu\cdot X\cong X\cdot Y_\mu$. Опять по уже доказанному мы можем заключить, что $\varPhi$ вполне ограничено в пространстве $Y_\mu\cdot X_\lambda$ при любом $\lambda\in\bf{Ord}$. И после этого, опять применяя теорему \ref{TH:X-cdot-Y-cong-Y-cdot-X}, мы можем сделать вывод, что $\varPhi$ вполне ограничено в пространстве $X_\lambda\cdot Y_\mu\cong Y_\mu\cdot X_\lambda$.
\epr

\blm\label{LM:o-drobyah}
Пусть $X$ и $Y$ псевдополны. Если $U\in {\mathcal U}(X_\lambda)$ и $V\in {\mathcal U}(Y_\mu)$ то $U\cdot V\in {\mathcal D}(X\cdot Y)$
\elm
\bpr
Пусть $\varPhi\subseteq X\cdot Y=Y\underset{\xi}{:}X^\star$ --- вполне ограниченное множество. По лемме \ref{LM:o-sohranenii-polnoi-ogranichennosti}, $\varPhi$ будет вполне ограничено и в пространстве $X_\lambda\cdot Y_\mu=Y_\mu\underset{\xi}{:}(X_\lambda)^\star$. Поскольку $U\cdot V=U:(V^\circ)$ --- окрестность нуля в пространстве $X_\lambda\cdot Y_\mu=Y_\mu\underset{\xi}{:}(X_\lambda)^\star$, найдется конечное множество $A\subseteq X_\lambda\cdot Y_\mu=Y_\mu\underset{\xi}{:}(X_\lambda)^\star$ такое что
\beq\label{o-drobyah-1}
\varPhi\subseteq U\cdot V+A.
\eeq
Но по лемме \ref{LM:o-drobyah-*}, все эти множества можно считать содержащимися в пространстве $X\cdot Y$. Мы получем такую вещь: для любого вполне ограниченного множества $\varPhi$ в $X\cdot Y$ существует конечное множество $A$ в $X\cdot Y$ такое, что выполняется \eqref{o-drobyah-1}. Это значит, что множество $U\cdot V$ является емким в пространстве $X\cdot Y$.
\epr

\blm\label{LM:(X_lambda-cdot-Y_lambda)^vartriangle=(X-cdot-Y)^vartriangle}
Для любых псевдополных локально выпуклых пространств $X$ и $Y$ и любых ординалов $\lambda,\mu\in\bf{Ord}$
\bit{
\item[(i)] при псевдонасыщении пространства $X_\lambda\cdot Y_\mu$ и $X\cdot Y$ становятся изоморфными:
\beq\label{(X_lambda-cdot-Y_mu)^vartriangle=(X-cdot-Y)^vartriangle}
(X_\lambda\cdot Y_\mu)^\vartriangle=(X\cdot Y)^\vartriangle,
\eeq

\item[(ii)] для любых функционалов $f\in (X_\lambda)^\star$ и $g\in (Y_\mu)^\star$ функционал $f\cdot g:X_\lambda\cdot Y_\mu\to\C$ непрерывен на пространстве $(X\cdot Y)^\vartriangle$, то есть  существует единственный функционал $h\in ((X\cdot Y)^\vartriangle)^\star$, замыкающий диаграмму
\beq\label{wedge_lambda^X-cdot-wedge_mu^Y}
\xymatrix @R=2.pc @C=4.0pc 
{
(X_\lambda\cdot Y_\mu)^\vartriangle \ar[d]_{\vartriangle_{X_\lambda\cdot Y_\mu}} \ar[r]^{(\wedge_\lambda^X\cdot\wedge_\mu^Y)^\vartriangle} & (X\cdot Y)^\vartriangle \ar@{-->}[d]^{h} \\
X_\lambda\cdot Y_\mu \ar[r]_{f\cdot g}& \C
}
\eeq
}\eit
\elm
\bpr
1. Здесь надо сразу заметить, что это утверждение достаточно доказать для случая $\mu=0$. Действительно, предположим, мы доказали утверждения
\bit{
\item[$(i)^\circ$] при псевдонасыщении пространства $X_\lambda\cdot Y$ и $X\cdot Y$ становятся изоморфными:
\beq\label{(X_lambda-cdot-Y)^vartriangle=(X-cdot-Y)^vartriangle}
(X_\lambda\cdot Y)^\vartriangle=(X\cdot Y)^\vartriangle,
\eeq

\item[$(ii)^\circ$] для любых функционалов $f\in (X_\lambda)^\star$ и $g\in Y^\star$ функционал $f\cdot g:X_\lambda\cdot Y\to\C$ непрерывен на пространстве $(X\cdot Y)^\vartriangle$, то есть существует единственный функционал $h\in ((X\cdot Y)^\vartriangle)^\star$, замыкающий диаграмму
\beq\label{wedge_lambda^X-cdot-id_Y}
\xymatrix @R=2.pc @C=4.0pc 
{
(X_\lambda\cdot Y)^\vartriangle \ar[d]_{\vartriangle_{X_\lambda\cdot Y}} \ar[r]^{(\wedge_\lambda^X\cdot\id_Y)^\vartriangle} & (X\cdot Y)^\vartriangle \ar@{-->}[d]^{h} \\
X_\lambda\cdot Y \ar[r]_{f\cdot g}& \C
}
\eeq
}\eit
Тогда, во-первых, будет выполняться \eqref{(X_lambda-cdot-Y_mu)^vartriangle=(X-cdot-Y)^vartriangle}, потому что
$$
(X_\lambda\cdot Y_\mu)^\vartriangle=\eqref{(X_lambda-cdot-Y)^vartriangle=(X-cdot-Y)^vartriangle}=(X\cdot Y_\mu)^\vartriangle=\eqref{X-cdot-Y-cong-Y-cdot-X}=
(Y_\mu\cdot X)^\vartriangle=\eqref{(X_lambda-cdot-Y)^vartriangle=(X-cdot-Y)^vartriangle}=
(Y\cdot X)^\vartriangle=\eqref{X-cdot-Y-cong-Y-cdot-X}=(X\cdot Y)^\vartriangle,
$$
и, во-вторых, для любых функционалов $f\in (X_\lambda)^\star$ и $g\in (Y_\mu)^\star$
функционал
$$
f\cdot g\ \circ\vartriangle_{X_\lambda\cdot Y_\mu}
$$
будет непрерывен на пространстве $(X_\lambda\cdot Y_\mu)^\vartriangle$,  а значит, и на изоморфном ему в силу только что доказанного равенства \eqref{(X_lambda-cdot-Y_mu)^vartriangle=(X-cdot-Y)^vartriangle} пространстве $(X\cdot Y)^\vartriangle$, и этот функционал на  $(X\cdot Y)^\vartriangle$ как раз и будет функционалом $h$, который строится в \eqref{wedge_lambda^X-cdot-wedge_mu^Y}. Его единственность также следует из \eqref{(X_lambda-cdot-Y_mu)^vartriangle=(X-cdot-Y)^vartriangle}.

2. Итак, мы поняли, что достаточно доказать более слабые утверждения $(i)^\circ$ и $(ii)^\circ$. Они доказываются индукцией по ординалам $\lambda\in\bf{Ord}$.

Прежде всего, при $\lambda=0$ равенство \eqref{(X_lambda-cdot-Y_mu)^vartriangle=(X-cdot-Y)^vartriangle} становится тривиальным
$$
(X\cdot Y)^\vartriangle=(X\cdot Y)^\vartriangle,
$$
а диаграмма \eqref{wedge_lambda^X-cdot-id_Y} превращается в диаграмму
$$
\xymatrix @R=2.pc @C=8.0pc 
{
(X\cdot Y)^\vartriangle \ar[d]_{\vartriangle_{X\cdot Y}} \ar[r]^{(\id_X\cdot\id_Y)^\vartriangle=\id_{X\cdot Y}^\vartriangle} & (X\cdot Y)^\vartriangle \ar@{-->}[d]^{h} \\
X\cdot Y \ar[r]_{f\cdot g}& \C
}
$$
в которой нужно просто положить $h=f\cdot g\circ\vartriangle_{X\cdot Y}$.

Далее, предположим, что мы доказали утверждения $(i)^\circ$ и $(ii)^\circ$ для всех ординалов $\iota<\lambda$, где $\lambda$ --- какой-то фиксированный ординал, то есть
\bit{
\item[$(i)^{\circ\circ}$] при псевдонасыщении пространства $X_\iota\cdot Y$ и $X\cdot Y$ становятся изоморфными:
\beq\label{(X_iota-cdot-Y)^vartriangle=(X-cdot-Y)^vartriangle}
(X_\iota\cdot Y)^\vartriangle=(X\cdot Y)^\vartriangle,
\eeq

\item[$(ii)^{\circ\circ}$] для любых функционалов $f\in (X_\iota)^\star$ и $g\in Y^\star$ функционал $f\cdot g:X_\iota\cdot Y\to\C$ непрерывен на пространстве $(X\cdot Y)^\vartriangle$, то есть существует единственный функционал $h\in ((X\cdot Y)^\vartriangle)^\star$, замыкающий диаграмму
\beq\label{wedge_iota^X-cdot-id_Y}
\xymatrix @R=2.pc @C=4.0pc 
{
(X_\iota\cdot Y)^\vartriangle \ar[d]_{\vartriangle_{X_\iota\cdot Y}} \ar[r]^{(\wedge_\iota^X\cdot\id_Y)^\vartriangle} & (X\cdot Y)^\vartriangle \ar@{-->}[d]^{h} \\
X_\iota\cdot Y\ar[r]_{f\cdot g}& \C
}
\eeq
}\eit
Покажем, что тогда $(i)^{\circ\circ}$ и $(ii)^{\circ\circ}$ выполняются для самого ординала $\lambda$ (подставленного вместо $\iota$). Здесь нужно рассмотреть два случая.

a) Предположим сначала, что $\lambda$ -- изолированный ординал, то есть
$$
\lambda=\iota+1
$$
для некоторого ординала $\iota<\lambda$. Покажем, что тогда для $\lambda$ выполняется условие $(ii)^\circ$. Пусть
$$
f\in (X_\lambda)^\star=(X_{\iota+1})^\star=((X_\iota)^\wedge)^\star=
\text{\cite[Theorem 3.10]{Ak03}}=
((X_\iota)^\star)^\vee
$$
и
$$
g\in Y^\star
$$
Тогда $f$ является пределом вполне ограниченной направленности $\{f_i;\ i\in I\}\subseteq (X_\iota)^\star$:
\beq\label{(X_lambda-cdot-Y_mu)^vartriangle=(X-cdot-Y)^vartriangle-1}
f_i\overset{((X_\iota)^\star)^\vee}{\underset{i\to\infty}{\longrightarrow}}f
\eeq
Поскольку множество $\{f_i\}\subseteq (X_\iota)^\star$ вполне ограничено, по следствию  \ref{LM:F-cdot-g-subseteq-(X-cdot-Y)^star}, множество
$$
\{f_i\cdot g;\ i\in I\}
$$
вполне ограничено в $(X_\iota\cdot Y)^\star$. С другой стороны, из теоремы \ref{TH:(f,g)|->f-cdot-g-razdelno-nepreryvno} следует, что $\{f_i\cdot g;\ i\in I\}$ является направленностью Коши. Итак, мы получаем, что $\{f_i\cdot g;\ i\in I\}$ --- вполне ограниченная направленность Коши в пространстве $(X_\iota\cdot Y)^\star$, значит, она сходится в объемлющем псевдополном пространстве $(X_\iota\cdot Y)^{\star\triangledown}$:
$$
(X_\iota\cdot Y)^\star\subseteq (X_\iota\cdot Y)^{\star\triangledown}=
\text{\cite[Theorem 3.14]{Ak03}}=(X_\iota\cdot Y)^{\vartriangle\star}=\eqref{(X_iota-cdot-Y)^vartriangle=(X-cdot-Y)^vartriangle}=
(X\cdot Y)^{\vartriangle\star}
$$
То есть найдется функционал $h\in (X\cdot Y)^{\vartriangle\star}$ такой что
\beq\label{(X_lambda-cdot-Y_mu)^vartriangle=(X-cdot-Y)^vartriangle-2}
f_i\cdot g\overset{(X\cdot Y)^{\vartriangle\star}}{\underset{i\to\infty}{\longrightarrow}}h
\eeq
С другой стороны, из цепочки равенств
$$
((X_\iota)^\star)^\vee=\text{\cite[Theorem 3.10]{Ak03}}= ((X_\iota)^\wedge)^\star=(X_{\iota+1})^\star=
(X_\lambda)^\star
$$
мы получаем, что соотношение \eqref{(X_lambda-cdot-Y_mu)^vartriangle=(X-cdot-Y)^vartriangle-1} эквивалентно соотношению
\beq\label{(X_lambda-cdot-Y_mu)^vartriangle=(X-cdot-Y)^vartriangle-3}
f_i\overset{(X_\lambda)^\star}{\underset{i\to\infty}{\longrightarrow}}f
\eeq
из которого по теореме \ref{TH:(f,g)|->f-cdot-g-razdelno-nepreryvno} следует
\beq\label{(X_lambda-cdot-Y_mu)^vartriangle=(X-cdot-Y)^vartriangle-4}
f_i\cdot g\overset{(X_\lambda\cdot Y)^{\star}}{\underset{i\to\infty}{\longrightarrow}} f\cdot g
\eeq
Вместе соотношения \eqref{(X_lambda-cdot-Y_mu)^vartriangle=(X-cdot-Y)^vartriangle-2} и \eqref{(X_lambda-cdot-Y_mu)^vartriangle=(X-cdot-Y)^vartriangle-4} означают, в частности, что на каждом операторе
$$
\ph\in X_\lambda\cdot Y=Y\underset{\xi}{:} (X_\iota)^\star=Y\underset{\xi}{:} X^\star=X_\lambda\cdot Y
$$
(по лемме \ref{LM:o-drobyah-*}, пространства $Y\underset{\xi}{:} (X_\iota)^\star$ и $Y\underset{\xi}{:} X^\star$ совпадают как множества, различаясь только топологией) функционалы $f\cdot g$ и $h$ совпадают:
$$
(f\cdot g)(\ph)=h(\ph).
$$
Это доказывает утверждение $(ii)^\circ$ (для случая $\lambda=\iota+1$).

Докажем теперь $(i)^\circ$. Отображение
$$
(X_\lambda\cdot Y)^\vartriangle=(Y\underset{\xi}{:}(X_\lambda)^\star)^\vartriangle\to (Y\underset{\xi}{:}X^\star)^\vartriangle=(X\cdot Y)^\vartriangle
$$
в прямую сторону всегда непрерывно. Нужно доказать, что оно непрерывно в обратную сторону, то есть, что непрерывно отображение
$$
(X\cdot Y)^\vartriangle=(Y\underset{\xi}{:}X^\star)^\vartriangle\to (Y\underset{\xi}{:}(X_\lambda)^\star)^\vartriangle=(X_\lambda\cdot Y)^\vartriangle
$$
Поскольку пространство $(X\cdot Y)^\vartriangle$ псевдонасыщено, нам достаточно доказать непрерывность отображения
$$
(X\cdot Y)^\vartriangle=(Y\underset{\xi}{:}X^\star)^\vartriangle\to Y\underset{\xi}{:}(X_\lambda)^\star=X_\lambda\cdot Y
$$
Рассмотрим базисную окрестность нуля в $X_\lambda\cdot Y$, то есть множество $U\cdot V$, где $U$ --- замкнутая выпуклая уравновешенная окрестность нуля в $X_\lambda$, а $V$ --- замкнутая выпуклая уравновешенная окрестность нуля в $Y$. По лемме \ref{LM:o-drobyah}, множество $U\cdot V$ является емким в $X\cdot Y$. Поэтому если мы докажем, что оно замкнуто в $(X\cdot Y)^\vartriangle$, то это будет означать, что оно является окрестностью нуля в $(X\cdot Y)^\vartriangle$.

Это становится очевидно, если представить $U\cdot V=V:U^\circ$ как поляру системы функционалов вида $\{f\cdot g;\ f\in U^\circ,\ g\in V^\circ\}$:
\beq\label{(X_lambda-cdot-Y_mu)^vartriangle=(X-cdot-Y)^vartriangle-5}
U\cdot V=V:U^\circ=\{f\cdot g;\ f\in U^\circ,\ g\in V^\circ\}^\circ
\eeq
Действительно,
\begin{multline*}
\ph\in U\cdot V=V:U^\circ\quad\Leftrightarrow\quad \ph(U^\circ)\subseteq V
\quad\Leftrightarrow\quad \forall f\in U^\circ \ \ph(f)\in V
\quad\Leftrightarrow \\ \Leftrightarrow\quad
\forall g\in V^\circ\ \forall f\in U^\circ \ |\underbrace{g(\ph(f))}_{\scriptsize \begin{matrix}\|\\ (f\cdot g)(\ph)\end{matrix}}|\le 1
\quad\Leftrightarrow\quad \ph\in \{f\cdot g;\ f\in U^\circ,\ g\in V^\circ\}^\circ
\end{multline*}
Теперь нужно заметить, что в формуле \eqref{(X_lambda-cdot-Y_mu)^vartriangle=(X-cdot-Y)^vartriangle-5}
$f\in U^\circ$, где $U$ --- окрестность нуля в $X_\lambda$, а $g\in V^\circ$, где $V$ --- окрестность нуля в $Y$. Значит, $f\in (X_\lambda)^\star$, а  $g\in Y^\star$. Отсюда по уже доказанному свойству $(ii)^\circ$ (для случая $\lambda=\iota+1$),  $f\cdot g$ --- непрерывные функционалы на пространстве $(X\cdot Y)^\vartriangle$. Мы получаем, что $U\cdot V$ является полярой некоторой системы непрерывных функционалов на $(X\cdot Y)^\vartriangle$, значит, $U\cdot V$ --- замкнутое множество в $(X\cdot Y)^\vartriangle$. Вдобавок, оно емкое, значит, оно является окрестностью нуля в (псевдонасыщенном) пространстве $(X\cdot Y)^\vartriangle$. Мы доказали свойство $(i)^\circ$ (для случая $\lambda=\iota+1$).

b) Теперь рассмотрим случай, когда $\lambda$ --- предельный ординал, то есть
$$
\lambda\ne \iota+1
$$
для всех $\iota<\lambda$. Здесь опять сначала нужно доказать $(ii)^\circ$. Пусть $f\in (X_\lambda)^\star$ и $g\in Y^\star$. Тогда
$$
f\in (X_\lambda)^\star=\text{\cite[Theorem 3.12]{Ak03}}
=(X^\star)^\lambda= \text{\cite[(1.32)]{Ak03}}=\bigcup_{\iota<\lambda}(X^\star)^\iota,
$$
поэтому существует ординал $\iota<\lambda$ такой, что $f\in (X^\star)^\iota=\text{\cite[Theorem 3.12]{Ak03}}=(X_\iota)^\star$. Мы получаем, что $f\in (X_\iota)^\star$ и $g\in Y^\star$, поэтому по предположению индукции $(ii)^{\circ\circ}$, функционал $f\cdot g$ должен быть непрывен на пространстве $(X\cdot Y)^\vartriangle$. Это доказывает  $(ii)^\circ$ (для случая предельного $\lambda$).

Теперь перейдем к $(i)^\circ$. Здесь повторяются рассуждения пункта a): нам нужно доказать непрерывность отображения
$$
(X\cdot Y)^\vartriangle=(Y\underset{\xi}{:}X^\star)^\vartriangle\to Y\underset{\xi}{:}(X_\lambda)^\star=X_\lambda\cdot Y
$$
Мы берем базисную окрестность нуля в $X_\lambda\cdot Y$, то есть множество $U\cdot V$, где $U$ --- окрестность нуля в $X_\lambda$, а $V$ --- окрестность нуля в $Y$, затем доказываем равенство
\eqref{(X_lambda-cdot-Y_mu)^vartriangle=(X-cdot-Y)^vartriangle-5}, затем замечаем, что в нем $f\in (X_\lambda)^\star$, а  $g\in Y^\star$. После этого по уже доказанному свойству $(ii)^\circ$ (для случая предельного $\lambda$), мы получаем, что $f\cdot g$ --- непрерывные функционалы на пространстве $(X\cdot Y)^\vartriangle$. Значит, $U\cdot V$ --- замкнутое множество в $(X\cdot Y)^\vartriangle$. Вдобавок, оно емкое, по лемме \ref{LM:o-drobyah}, значит, оно является окрестностью нуля в (псевдонасыщенном) пространстве $(X\cdot Y)^\vartriangle$. Мы доказали свойство $(i)^\circ$ (для случая предельного $\lambda$).
\epr

Для произвольных локально выпуклых пространств $X$ и $Y$ отображения
$$
\vartriangle_X: X^\vartriangle\to X,\qquad \vartriangle_Y: Y^\vartriangle\to Y
$$
определяют отображение
\beq\label{vartriangle_X-cdot-vartriangle_Y}
(\vartriangle_X\cdot \vartriangle_Y): (X^\vartriangle\cdot Y^\vartriangle) \to (X\cdot Y)
\eeq
(которое просто совпадает с отображением \eqref{wedge_lambda^X-cdot-wedge_mu^Y-1} при подходящих $\lambda$ и $\mu$).

\btm\label{TH:(X^vartriangle-cdot-Y^vartriangle)^vartriangle=(X-cdot-Y)^-vartriangle}
Для любых псевдополных локально выпуклых пространств $X$ и $Y$ псевдонасыщение отображения \eqref{vartriangle_X-cdot-vartriangle_Y} является изоморфизмом локально выпуклых пространств:
\beq\label{(X^vartriangle-cdot-Y^vartriangle)^vartriangle=(X-cdot-Y)^-vartriangle}
X^\vartriangle\odot Y^\vartriangle\cong (X^\vartriangle\cdot Y^\vartriangle)^\vartriangle\cong (X\cdot Y)^\vartriangle.
\eeq
\etm
\bpr
Здесь второй изоморфизм получается из \eqref{(X_lambda-cdot-Y_mu)^vartriangle=(X-cdot-Y)^vartriangle}, если подобрать ординалы $\lambda$ и $\mu$ так, чтобы $X_\lambda=X^\vartriangle$ и $Y_\mu=Y^\vartriangle$. А в первом изоморфизме используется тот факт, что псевдонасыщение $X^\vartriangle$ псевдополного пространства $X$ всегда псевдополно (и поэтому стереотипно) \cite[Proposition 3.16]{Ak03}. Отсюда следует, что в пространстве $(X^\vartriangle)^\star$ всякое вполне ограниченное множество $S$ содержится в поляре $U^\circ$ некоторой окрестности нуля $U\subseteq X^\vartriangle$. Как следствие, в пространстве операторов $Y^\vartriangle:(X^\vartriangle)^\star$ топология равномерной сходимости на вполне ограниченных множествах совпадает с топологией равномерной сходимости на полярах $U^\circ$ окрестностей нуля $U\subseteq X^\vartriangle$:
$$
Y^\vartriangle:(X^\vartriangle)^\star=Y^\vartriangle\underset{\xi}{:}(X^\vartriangle)^\star
$$
Отсюда мы получаем цепочку, доказывающую первый изоморфизм в \eqref{(X^vartriangle-cdot-Y^vartriangle)^vartriangle=(X-cdot-Y)^-vartriangle}:
$$
X^\vartriangle\odot Y^\vartriangle=Y^\vartriangle\oslash (X^\vartriangle)^\star=
(Y^\vartriangle:(X^\vartriangle)^\star)^\vartriangle=
(Y^\vartriangle\underset{\xi}{:}(X^\vartriangle)^\star)^\vartriangle= \eqref{DEF:X-cdot-Y}=(X^\vartriangle\cdot Y^\vartriangle)^\vartriangle.
$$
\epr

\tableofcontents